\documentclass[12pt]{article}

\usepackage[latin1]{inputenc}

   \usepackage[francais,english]{babel}

   \usepackage{graphicx} 

\usepackage{amsmath}    
\usepackage{amssymb}

\input diagrams.sty
\usepackage[all]{xy} 

\begin{document}

\title{\bf{Singularit\'es r\'eelles isol\'ees et d\'eveloppements asymptotiques d'int\'egrales oscillantes}}

\author{Daniel Barlet}

\date{Universit\'e Nancy I  et Institut Universitaire de France}

\maketitle

\noindent \textit{ De la fen\^etre de mon bureau je contemple le vieux pont de Malz\'eville et, me projetant  65  ans en arri\`ere, j'imagine Jean Leray, alors jeune professeur \`a l'Universit\'e de Nancy, regagnant son domicile depuis la porte de le Craffe, pr\`es de laquelle \'etait alors l'Institut de Math\'ematique, s'arr\^etant sur ce vieux pont pour scruter les tourbillons de la Meurthe derri\`ere les piles du pont, et s'interro\-geant sur les solutions turbulentes  des \'equations de Navier-Stokes .}

\bigskip

\bigskip

\noindent{\bf Abstract}

\noindent Let  $ (X_{\mathbb{R}} ,0) $ \ be a germ of  real analytic subset in  $ (\mathbb{R}^{N} ,0) $  of pure dimension $n+1 $  with an isolated singularity at  $0$ . Let  
 $$ (f_{\mathbb{R}},0) :  (X_{\mathbb{R}} ,0) \longrightarrow  (\mathbb{R},0) $$
a real analytic germ with an isolated singularity  at  $0$ , such that its complexification  $f_{\mathbb{C}}$ \ vanishes on the singular set  $ S $  of  $ X_{\mathbb{C}} $ .We also assume that  $X_{\mathbb{R}} - \lbrace 0 \rbrace $  \ is orientable .

\noindent To each \   $  A \in H^{0}(X_{\mathbb{R}} - \lbrace 0 \rbrace ,\mathbb{C} ) $ \ we associate  a $n-$cycle  $ \Gamma(A) $  ("explicitly " described) in the complex  Milnor fiber of  $f_{\mathbb{C}}$ \ at  $0$  such that the non trivial terms in the asymptotic expansions of the oscillating integrals \  $ \int_{A} e^{i\tau f(x)} \varphi(x) $ \  when  $ \tau \rightarrow  \pm \infty $ \ can be read from the spectral decomposition of  \ $\Gamma(A) $ \ relative to the monodromy of  \ $f_{\mathbb{C}}$ \ at  \ $0$ .

\bigskip

\newpage

\noindent{\bf R\'esum\'e}

\bigskip

\noindent Soit $ (X_{\mathbb{R}} ,0) $ \ un germe de sous-ensemble analytique r\'eel \`a l'origine de  $ \mathbb{R}^{N}  $ de dimension pure  $n+1$ \ ayant une singularit\'e isol\'ee en $0$ . Soit 
 $$ (f_{\mathbb{R}},0) :  (X_{\mathbb{R}} ,0) \longrightarrow  (\mathbb{R},0) $$
un germe de fonction analytique r\'eelle ayant une singularit\'e isol\'ee en  $0$  telle que sa complexifi\'ee  $f_{\mathbb{C}}$ \ s'annule sur le lieu singulier   $ S $ de  $ X_{\mathbb{C}} $ . Nous supposerons \'egalement que la vari\'et\'e analytique r\'eelle $X_{\mathbb{R}} - \lbrace 0 \rbrace $\ est  orientable.

\noindent A chaque   $  A \in H^{0}(X_{\mathbb{R}} - \lbrace 0 \rbrace ,\mathbb{C} ) $ \  nous associons un  $n-$cycle  $ \Gamma(A) $  (explicitement d\'ecrit) dans la fibre de Milnor complexe de  $f_{\mathbb{C}}$ \    
 en  $0$  tel que les termes non triviaux dans les d\'eveloppements asymptotiques  quand $ \tau \rightarrow  \pm \infty $ \  des int\'egrales oscillantes  \  $ \int_{A} e^{i\tau f(x)} \varphi(x) $ \ soient d\'etect\'es par la d\'ecomposition spectrale de  $\Gamma(A) $ \  par rapport \`a la monodromie de  \ $f_{\mathbb{C}}$ \ en  \ $0$ .

\bigskip

\bigskip

\bigskip

\noindent{\bf Table des mati\`eres}

\begin{itemize}

\item{\bf 1. Introduction.}

\item{\bf 2. Transformation de Mellin sur  $ \mathbb{R}^{*}. $}

\item {\bf 3. Cohomologie relative et variation .}

\item{\bf 4. Le cas d'une valeur propre $ \not = 1 $ .}

\item{\bf 5. Le cas de la valeur propre  1 .}

\item{\bf 6. Le cas o\`u  $ \partial A \subset \lbrack 0 \rbrack .$ }

\item{\bf 7. Annexe .}

\item{\bf R\'ef\'erences .}

\end{itemize}

\newpage

\section{\bf Introduction}

Ce texte est un hommage \`a Jean Leray ; on constatera facilement que par les id\'ees et les techniques utilis\'ees, il est un des nombreux descendants en ligne directe du fameux article  [L.59]  "le probl\`eme de Cauchy III" qui est l'une des grandes contributions de Jean Leray \`a la th\'eorie des fonctions de plusieurs variables complexes.
Cet article peut \^etre lu  \`a deux niveaux . Le premier niveau, qui est proche de la conf\'erence donn\'ee \`a Nantes, est une introduction \`a la th\'eorie des singularit\'es  d'une fonction analytique r\'eelle ; son but est d'expliquer comment celle-ci permet de d\'ecrire les d\'eveloppements asymptotiques des int\'egrales oscillantes \`a phase analytique r\'eelle. Il est assez surprenant et int\'eressant de voir que la topologie de l'application d\'efinie par la complexifi\'ee de la phase consid\'er\'ee et  la topologie de la position du r\'eel dans le complexe  suffisent \`a pr\'evoir exactement quels types de termes appara\^itront dans ces d\'eveloppements asymptotiques. Dans cette optique, on peut simplement consid\'erer le cas o\`u $ ({X }_{\mathbb{R}},0) =  (\mathbb{R}^{n+1 },0)$  et se contenter de lire seulement les \'enonc\'es des r\'esultats ainsi que les constructions qui les pr\'ec\`edent d\'efinissant les cycles  $ \Gamma(A) $ \ et \ $\widehat{\Gamma}(A) $ \ respectivement pour les th\'eor\`emes 1, 1bis  et  2 . Les \'enonc\'es des corollaires \'etant probablement plus faciles \`a comprendre que ceux  des th\'eor\`emes eux-m\^eme .
L'autre niveau est celui d'un article de recherche qui  am\'eliore les r\'esultats d\'ej\`a connus sur ce sujet .En particulier nous g\'en\'eralisons substantiellement les r\'esultats de [J.91] ,[B.M.02 ]  et [B.02]  . 

\noindent De fa{\c c}on pr\'ecise, nous d\'ecrivons , si $ ({X }_{\mathbb{R}},0)\subset  (\mathbb{R}^{N},0)$ est un germe d'ensemble analytique r\'eel de dimension pure $( n+1)$  \`a singularit\'e isol\'ee en 0 (donc nous \'etudions des probl\`emes de "phase stationnaire" avec une singularit\'e ambiante l\`a o\`u la phase stationne) et si $ (f_{\mathbb{R}},0) :  ({X }_{\mathbb{R}},0) \longrightarrow ((\mathbb{R}, 0)$  est un germe analytique r\'eel \`a singularit\'e isol\'ee en 0 sur $ ({X }_{\mathbb{R}},0)$ , sous l'hypoth\`ese que le lieu singulier du complexifi\'e $ ({X }_{\mathbb{C}},0)  $ est contenu dans $\lbrace f_{\mathbb{C}}=0 \rbrace $ et que ${X }_{\mathbb{R}}-\lbrace 0 \rbrace  $ est orientable, les termes qui vont appara\^itre dans le d\'eveloppement asymptotique quand $ \tau \longrightarrow + \infty $ de l'int\'egrale 
$$\int_{A}^{ }e^{i\tau f(x)}.\varphi (x) $$
 o\`u $\varphi $ est une forme \ $\mathcal{C}^{\infty }$ \ de degr\'e $ n+1$ \ ,  \`a support compact (une forme- test) , et o\`u  $ {A}=\sum_{\alpha} a_{\alpha} .A_{\alpha} \in  H^{0}(X_{\mathbb{R}}-f_{\mathbb{R}}^{-1}(0),\mathbb{C})$ est une combinaison lin\'eaire \`a coefficients complexes de composantes connexes de $X_{\mathbb{R}}-f_{\mathbb{R}}^{-1}(0) $ . Nous exhibons un n-cycle $ \Gamma(A) $ associ\'e \`a  A  dans la fibre de Milnor  en 0 de la complexifi\'ee $f_{\mathbb{C}}$ de $f_{\mathbb{R}}$ dont la d\'ecomposition spectrale par rapport \`a la monodromie en 0 de $f_{\mathbb{C}}$ gouverne ces d\'eveloppements asymptotiques. 

\noindent Les \'enonc\'es pr\'ecis sont donn\'es aux th\'eor\`emes 1, 1bis et  2  et dans les corollaires qui les suivent . Evidemment, si l'on sait montrer que pour un $ A $ \  donn\'e le cycle $ \Gamma(A) $  n'est pas nul dans la cohomologie de la fibre de Milnor, on en d\'eduit qu'il existe une forme test $\varphi $ pour laquelle la fonction  $\tau \longrightarrow \int_{A}^{} e^{i\tau f(x)}.\varphi (x) $ n'est pas une fonction de la classe de L. Schwartz sur $ \mathbb{R} $. C'est en fait la m\'ethode utilis\'ee par A. Jeddi dans  [J.02] (en s'inspirant de l'article "fondateur" [M.74])  pour r\'esoudre la "conjecture" de Palamodov (voir [P.86]), c'est \`a dire pour montrer que dans le cas o\`u  $ ({X }_{\mathbb{R}},0)=(\mathbb{R}^{ n+1},0)$ et o\`u ${A}=\sum_{\alpha} A_{\alpha} $ (on int\`egre donc sur $\mathbb{R}^{n+1}$) il existe une telle forme  test  $\varphi .$ En fait , A. Jeddi d\'eduit le cas g\'en\'eral d'une fonction analytique r\'eelle sur $\mathbb{R}^{ n+1}$ du cas d'une fonction \`a singularit\'e isol\'ee.

 \noindent  Pour terminer cette introduction, je voudrais remercier le Laboratoire de Math\'ematiques de l'Universit\'e de Nantes, maintenant Laboratoire Jean Leray, pour m'avoir invit\'e \`a donner cette conf\'erence lors de son "bapt\^eme" .

\bigskip

\section{\bf Transformation de Mellin sur $ \mathbb{ R}^{*} $ }

\noindent   Soit  $\varphi $ une fonction  $\mathcal{C}^{\infty }$ sur $\mathbb{R}^* $ qui est born\'ee et \`a support born\'e ; nous d\'efinirons la transform\'ee de Mellin de $\varphi $ pour $\Re(\lambda) > 0 $  \ par la formule suivante :  
 
 $$ M\varphi  (\lambda) := \frac{1}{i\pi}[\int_{0}^{+\infty} x^{\lambda}\varphi (x) \frac{dx}{x}  - e^{-i\pi \lambda}.\int_{0}^{+\infty} x^{\lambda}\varphi (-x) \frac{dx}{x}] $$

\bigskip

 \noindent C'est un exercice \'el\'ementaire de voir que si  $\varphi $  se prolonge en une fonction  $\mathcal{C}^{\infty }$ sur   $\mathbb{R}$, alors  $ M\varphi $ se prolonge en une fonction enti\`ere (voir [B.99] ).

\noindent Pour illustrer cette d\'efinition, donnons un r\'esultat qui nous sera utile plus loin :

\bigskip

\noindent {\bf Lemme 1.}

\noindent {\it  Soit  $ s_0 >0 $  et   soit  $\varphi $ la fonction d\'efinie par :
\bigskip

\noindent   \    $\varphi (s) = (s/s_{0})^{r}.P[Log(s/s_{0})]$ \  \  \  \    \  pour  \  \ $  0< s < s_0 $ 

\bigskip

\noindent  \    $\varphi (s) = (s/s_{0})^{r}.Q[Log(s/s_{0})]$ \  \  \  \ pour  \  \ $  -s_0 < s < 0 $
\bigskip

\noindent  o\`u P et Q  sont des polyn\^omes et o\`u la partie imaginaire du logarithme est dans
l'intervalle \   $]-3\pi  /2,\pi /2 [ $ pour  z  n'appartenant pas \`a  $ i\mathbb{R}^{+ } $ . Alors le r\'esidu
en   $ \lambda   = -r $ de la fonction m\'eromorphe  :

 \noindent  $$ F (  \lambda) =  \int_{0}^{s_0} (s/s_{0})^{\lambda} \ \varphi (s) \frac{ds}{s}  - e^{-i\pi \lambda}.\int_{0}^{s_0} (s/s_{0})^{\lambda} \ \varphi (-s) \frac{ds}{s} $$

\noindent  est \'egal  \`a  P(0) - Q(0) . } 

\noindent {\bf Preuve.}

\noindent  Posons  \  $x =  s/s_{0} $ \ . Cela donne :

\noindent  $$ F(\lambda) =  \int_{0}^{1} x^{\lambda}x^{r}.P[Log x ] \frac{dx}{x}  - e^{-i\pi  \lambda}.\int_{0}^{1} x^{\lambda}.e^{-i\pi r}x^{r}.Q[Log x  - i\pi ]  \frac{dx}{x} $$

\noindent  Mais si  $ \Gamma $ d\'esigne le contour form\'e du segment  [-1,1] et du demi-cercle unit\'e inf\'erieur parcouru dans le sens indirect la nullit\'e de l'int\'egrale  :

$$  \int_{\Gamma}^{}  z^{\lambda + r}.Q[Log z ]\frac{dz}{z }  $$

\noindent et le fait que la contribution du demi-cercle \`a cette int\'egrale est une fonction enti\`ere de $\lambda$ montrent que le r\'esidu en $\lambda = - r $ de  $F(\lambda)$  est le m\^eme que celui de la fonction m\'eromorphe :

$$  \int_{0}^{1} x^{\lambda + r}.(P - Q)[Log x ]. \frac{dx}{x}   \quad .$$

\noindent On conclut alors facilement  \ \ \ \ \ \ \quad \qquad  $ \blacksquare$

\smallskip

\smallskip

\smallskip
\vskip 1cm

\noindent {\bf Remarque.}

\noindent  Ce lemme met en \'evidence une erreur de calcul dans la preuve du th\'eor\`eme 6.1 de [B.M.02] qui conduit aux formules erronn\'ees (6.2a) et (6.2b) de cet article. Les formules correctes sont donn\'ees aux th\'eor\`emes 1 et 1bis (pour plus de pr\'ecisions se reporter \`a la "Remarque/Erratum" qui suit le corollaire du th\'eor\`eme 1 ).

\bigskip


\noindent Consid\'erons maintenant la situation d'un germe  :

\noindent $(f_{\mathbb{R}},0) :  ({X }_{\mathbb{R}},0) \longrightarrow (\mathbb{R}, 0)$

\noindent de fonction analytique r\'eelle sur un espace analytique r\'eel  $ ({X}_{\mathbb{R}},0) $ satisfaisant les hypoth\`eses suivantes :

\begin{itemize}
\item [\it H a)]  ${(X}_{\mathbb{R}},0) $  est de dimension pure  $ n+1 $  et est lisse en dehors de l'origine .

\item [\it H b)]  $f_{\mathbb{R}}$ \ a une singularit\'e isol\'ee en  $ 0 $  sur  $ {X}_{\mathbb{R}} $

\item [\it H c)]  Le complexifi\'e  ${X}_{ \mathbb{C}} $ de ${ X}_{\mathbb{R}} $  a un lieu singulier  $ S $  contenu
                      dans  $ \lbrace f_{\mathbb{C} } =  0  \rbrace $  o\`u  $(f_{\mathbb{C}},0) :  ({X }_{\mathbb{C}},0) \longrightarrow (\mathbb{C}, 0)$  est la complexifi\'ee de $(f_{\mathbb{R}},0) $.

\item[\it H d)]  Notre derni\`ere hypoth\`ese sera que  ${X}_{\mathbb{R}}  - \lbrace  0 \rbrace $   est orientable , et nous fixerons une fois pour toute une orientation de cette vari\'et\'e analytique r\'eelle lisse de dimension $ n+1 $ .
\end{itemize}

\noindent Le lecteur non familier avec les espaces analytiques singuliers pourra supposer que l'on a simplement  $ {X}_{\mathbb{R}} = \mathbb{R}^{n+1} $ . Les r\'esultats sont d\'ej\`a int\'eressants dans ce cas (et les hypoth\`eses  H  a) , H c)  et  H d)  sont alors  automatiquement satisfaites ) .

\noindent Pr\'ecisons que, contrairement \`a la situation qui est  consid\'er\'ee dans  [B.M.02]  et  [B.02] , nous ne supposons pas ici que  le complexifi\'e  ${X}_{ \mathbb{C}} $  de  $ {X}_{\mathbb{R}}$ a une singularit\'e isol\'ee en  $ 0 $ . Donc  $ S $  peut \^etre de dimension positive, mais ne rencontre  ${ X}_{\mathbb{R}} $  qu'au plus \`a  l'origine . Nous ne supposons pas non plus que  $f_{\mathbb{C}}$  a une singularit\'e isol\'ee en $ 0 $  dans  ${X}_{ \mathbb{C}} $  (ce qui, pour  $ S \not=  \lbrace  0 \rbrace $  , n'a d'ailleurs pas un sens \'evident, mais qui est l'hypoth\`ese faite"en plus "de   $ S =  \lbrace  0 \rbrace $   dans   [B.M.02]  et  [B.02] ) .

\noindent Fixons maintenant un repr\'esentant de Milnor $$  f :  {X}_{\mathbb{C}}\longrightarrow  {D}  $$du germe  $(f_{\mathbb{C}},0)$  dans   ${(X}_{\mathbb{C}},0) $ (voir au  paragraphe 3 
 pour plus de pr\'ecisions sur cela ). Notons par   $ {X}_{\mathbb{R}}$  et   $f_{\mathbb{R}}$ les  repr\'esentants correspondants des germes  $({X}_{\mathbb{R}},0) $ et  $(f_{\mathbb{R}},0) $ . Soit  $ {A} = \sum_{\alpha } a_{\alpha}. A_{\alpha} $ un  \'el\'ement de  $H^{0}({X}_{\mathbb{R}} - f_{\mathbb{R}}^{-1}(0) , \mathbb{C} ) $ , c'est -\`a-dire une combinaison lin\'eaire \`a coefficients complexes de composantes connexes de  ${X}_{\mathbb{R}} - f_{\mathbb{R}}^{-1}(0) $  (il n'y a qu'un nombre fini de telles composantes connexes comme on le verra plus loin) . D\'efinissons  alors le 1-courant sur  ${ X}_{\mathbb{R}}  ,$ d\'ependant holomorphiquement du param\`etre  $ \lambda $ \ , en posant pour  $\Re{\lambda} \gg 1 $  
$$ \Box \longrightarrow \frac{1}{i\pi}  \int_{A}^{}  f^{\lambda} \  \Box \wedge \frac{df}{f}  $$
o\`u  $ \Box $  d\'esigne une $ n-$ forme  $\mathcal{C}^{\infty} $ \`a support compact dans  ${X}_{\mathbb{R}} $ ; c'est -\`a-dire que pour chaque  $ \Box $ donn\'ee, la valeur sur cette forme-test de ce 1-courant est  la transform\'ee de Mellin  (d\'efinie ci-dessus) de la fonction sur  $ \mathbb{R}^{*} $  d\'efinie par  :

\newpage

$$  s \longrightarrow \int_{{A}\bigcap \lbrace f_{\mathbb{R}} = s \rbrace }^{} \Box \qquad \qquad . $$

\noindent L'int\'egrale oscillante  $  \int_{A}^{} e^{i\tau f} \ \Box  \wedge df $ \  sera alors, par d\'efinition , la transform\'ee de Fourier de cette m\^eme fonction sur  $ \mathbb{R} $  (elle  est localement born\'ee \`a l'origine et \`a support compact ) .

\noindent Il est bien connu que, gr\^ace au th\'eor\`eme de d\'esingularisation  [H.64] de H. Hironaka,  la fonctiond\'efinie ci-dessus est   $\mathcal{C}^{\infty} $   sur  $ \mathbb{R}^{*} $ \  et admet  quand  $  s \rightarrow  \pm \infty $  des d\'eveloppements asymptotiques ind\'efiniment d\'erivables dans l'\'echelle des fonctions  $  s^{\alpha+ \nu}.(Log |s| )^{j}  $ \  avec \ $ \alpha \in [0,1[ \bigcap \mathbb{Q} $ \ , \ $ j \in [0,n] \bigcap \mathbb{N}  $  \ et  \ $  \nu \in \mathbb{N} $  (voir  [A.70] \  et \   [JQ.70] ).

\noindent On en d\'eduit imm\'ediatement que la distribution  $ \int_{A}^{}  f^{\lambda} \  \Box   $ ainsi que le 1-courant  $ \int_{A}^{}  f^{\lambda} \  \Box \wedge \frac{df}{f}  $ admettent des prolongements m\'eromorphe en $ \lambda $  \`a tout le plan complexe, avec un nombre fini de s\'eries de p\^oles d'ordre  $ \leq n+1 $  en  $ - \alpha - \nu $  o\`u  $ \alpha $  prend un nombre fini de valeurs dans  $ [0,1[ \bigcap \mathbb{Q}  $ et o\`u $ \nu \in \mathbb{N} $. De fa{\c c}on \'equivalente, on en d\'eduit que l'int\'egrale oscillante $  \int_{A}^{} e^{i\tau f} \ \Box  $   admet quand  $ \tau  \rightarrow  \pm \infty $ un d\'eveloppement asymptotique dans l'\'echelle des fonctions  $ \tau ^{- \alpha - \nu} (Log | \tau | )^{j} $ o\`u  $ \alpha ,j  $  et \ $ \nu  $  prennent les m\^emes valeurs  que plus haut  .

\noindent Nous formulerons nos r\'esultats en utilisant les p\^oles du prolongement m\'ero\-morphes de la distribution  $ \int_{A}^{} f^{\lambda} \  \Box   $  mais un dictionnaire facile
 (voir par exemple  [B.M.93]  prop.(0.7) ) 
 permet la traduction en terme des d\'eveloppements asymptotiques de l'int\'egrale oscillante  $  \int_{A}^{} e^{i\tau f} \ \Box  \wedge df $  , ces deux formulations \'etant \'equivalentes \`a l'\'etude des d\'eveloppements asymptotiques  \`a l'origine de { \bf "l'int\'egrale-fibre" } o\`u  \  $ \Box \in \mathcal{C}^{\infty}_{c}(X_{\mathbb{R}})^{n}  : $
 $$  s \longrightarrow  \int_{{A}\bigcap \lbrace f_{\mathbb{R}} = s \rbrace } \Box   \qquad \qquad .$$

\newpage

\section {\bf Cohomologie relative  et  variation .}

\bigskip

\noindent Soit  $( {X}_{\mathbb{C}},0) $  un germe d'ensemble analytique complexe irr\'eductible de dimension  n+1  dans 
 $ ( \mathbb{C}^{N}, 0 ) $ et soit 
 $ f =  f_{ \mathbb{C}} : ( \mathbb{C}^{N} , 0 )  \longrightarrow   ( \mathbb{C} , 0 ) $
 un germe de fonction holomorphe, nul sur le lieu singulier  de  $( {X}_{\mathbb{C}},0) $ \ , mais non identiquement nul sur \ $( {X}_{\mathbb{C}},0) $ . On notera par  $ \widehat{S} $  la r\'eunion du lieu singulier \ $ S $ \ de  $( {X}_{\mathbb{C}},0) $  et des points lisses de  $( {X}_{\mathbb{C}},0) $  en lesquelles la diff\'erentielle de  f  est nulle . On a  $\widehat{S }\subset  f^{-1}(0) $ sous ces hypoth\`eses  et pour \  $  0 < \varepsilon \ll 1 , 0 < \delta \ll \varepsilon $  \ on d\'efinit un repr\'esentant de Milnor en posant  :
 $${X}_{\mathbb{C}} : = \widetilde{ {X}_{\mathbb{C}}} \cap B(0, \varepsilon) \cap f^{-1}(D_{\delta})  $$  
\noindent o\`u $ \widetilde{ {X}_{\mathbb{C}}} $ d\'esigne un repr\'esentant du germe  $( {X}_{\mathbb{C}},0) $. La restriction de \  f \  \`a l'ouvert  $ {X}_{\mathbb{C}} - f^{-1}(0) $ est alors une fibration  $ \mathcal{C}^{\infty} $  sur le disque point\'e  $ D_{\delta}^{*} $  \  de fibre  $ F := f^{-1}(s_{0}) $ \  o\`u \  $ s_{0} \in D_{\delta} \cap \mathbb{R}^{+*} $  est un point base de $ D_{\delta}^{*} $ .

\noindent On remarquera que dans la construction de Milnor (voir [Mi.68] ) le champ de vecteur  $ \mathcal{C}^{\infty} $  peut \^etre prolong\'e de fa{\c c}on  $ \mathcal{C}^{\infty} $  au voisinage d'un compact  $ \Lambda $  de  $ f^{-1}(0) $  pourvu que l'on ait  $ \Lambda \cap \widehat{S} = \emptyset $ .

\noindent Un tel compact   $ \Lambda $  \'etant fix\'e ( seul le cas o\`u   $ \Lambda $  est une sous-vari\'et\'e analytique r\'eelle lisse de  $ f^{-1}(0) $  v\'erifiant  $ \Lambda \cap \widehat{S} = \emptyset $ nous int\'eressera ici ) on peut trouver des voisinages ouverts  $ \mathcal{U} \subset \mathcal{U}' $  de   $ \Lambda $ , d'adh\'erences disjointes de \ $ \widehat{S} $ , et v\'erifiant les propri\'et\'es suivantes  (quitte a restreindre $ \delta $ ) :

\begin{itemize}
\item[1)]   $ f : \mathcal{U}'  \longrightarrow   D_{\delta} $ \   est une fibration  $ \mathcal{C}^{\infty} $  triviale de fibre  $  F \cap \mathcal{U}'  $  et on a une trivialisation  $ \mathcal{C}^{\infty} $ \ 
$$  \Phi ' :   \mathcal{U}'  \longrightarrow   ( F \cap \mathcal{U}'  ) \times D_{\delta} $$  de cette fibration  qui induit un diff\'eomorphisme not\'e \ $\Phi$ \ de  \  $ \mathcal{U} $  sur   \ $   ( F \cap \mathcal{U} ) \times D_{\delta} $ .

\item[2)]  l'inclusion   $ F \cap \mathcal{U} \subset  F \cap \mathcal{U}'  $  est une \'equivalence d'homotopie  ( en fait dans le cas o\`u  $ \Lambda $  est une sous-vari\'et\'e compacte lisse de  $ f^{-1}(0)  -  \widehat{S} $ on peut choisir  $ F \cap \mathcal{U} $ \  et   \  $ F \cap \mathcal{U}' $ \ homotopiquement \'equivalents \`a  $ \Lambda $ , $ \mathcal{U} $  et   $ \mathcal{U}' $ \'etant les traces sur \ $ f^{-1}(D_{\delta}) , \delta $ assez petit, de voisinages tubulaires   $ \mathcal{C}^{\infty} $  de  $ \Lambda $ ) .

\item[3)] On a une trivialisation \ $\Psi $ \ de classe \   $ \mathcal{C}^{\infty} $  \ de la fibration de Milnor au dessus de l'ouvert simplement connexe  $ D_{\delta} - i \mathbb{R}^{+} \cap D_{\delta} $  qui est compatible \`a  $ \Phi' $ , c'est-\`a-dire que  $ \Phi' $ \  et  \ $ \Psi $  coincident sur l'ouvert   $$ \mathcal{U}' - f^{-1}(i \mathbb{R}^{+} \cap D_{\delta} ) \quad . $$

\end{itemize}

\noindent {\bf Cons\'equences importantes :}

\begin{itemize}

\item[a)]   \   \ Dans cette situation, la monodromie de  \ f \  est l'identit\'e sur   $ \mathcal{U}' $  et sur  $ \mathcal{U} $  . Elle agit donc comme l'identit\'e sur les cohomologies de  $ F \cap \mathcal{U}' $ \ et \  $ F \cap \mathcal{U} $  respectivement.
La monodromie agit \'egalement sur les cohomologies relatives (voir l'Annexe pour des pr\'ecisions sur les cohomologies de de Rham relatives \`a un ouvert que nous utiliserons ici )  $ H^{*}(F,F \cap \mathcal{U}) $  \  et  \ $  H_{c}^{*}(F,F \cap \mathcal{U})  .$

\item[b)]  On a des suites exactes d'espaces vectoriels monodromiques (voir  \'egale\-ment l'Annexe )
$$ \cdots  \rightarrow  H^{n-1}(F\cap \mathcal{U})   \rightarrow  H^{n}(F,F\cap \mathcal{U})  \rightarrow  H^{n}(F)  \rightarrow  H^{n}(F\cap \mathcal{U}) \rightarrow  \cdots $$
$$\cdots \rightarrow  H_{c}^{n}(F\cap \mathcal{U})  \rightarrow  H_{c}^{n}(F)   \rightarrow  H_{c}^{n}(F,F\cap \mathcal{U})  \rightarrow  H_{c}^{n+1}(F\cap \mathcal{U})  \rightarrow  \cdots $$

\noindent Si  $ V_{\not=1} $  d\'esigne la somme des sous-espaces spectraux pour les valeurs propres  $ \not=1 $  de la monodromie (resp.  $  V_{=1} $   le sous-espace spectral pour la valeur propre   1  ) agissant sur l'espace vectoriel monodromique  $ V $ \  l'application canonique :
$$  can :  H^{n}(F,F\cap \mathcal{U}) _{\not=1} \rightarrow  H^{n}(F) _{\not=1} $$
induit un isomorphisme puisque la monodromie agit comme l'identit\'e sur la cohomologie de  $ F\cap \mathcal{U} $. Il en est de m\^eme pour l'application 
$$can_{c}  :  H_{c}^{n}(F) _{\not=1}  \rightarrow  H_{c}^{n}(F,F\cap \mathcal{U}) _{\not=1}  \ \qquad \qquad  \qquad \Box $$

\end{itemize} 
\noindent Pour la valeur propre 1  nous allons construire une application de variation 
$$  var_{c} :  H_{c}^{n}(F,F\cap \mathcal{U}) _{=1}  \rightarrow  H_{c}^{n}(F)_{=1} $$
qui v\'erifiera les relations  : $$ var_{c} \circ can_{c} = T_{c} - 1  \  , \  can_{c} \circ var_{c} = T_{c} - 1 .$$

\noindent Nous allons utiliser pour cela la description de  $ H_{c}^{n}(F,F\cap \mathcal{U}) $ par des n-cha\^ines compactes orient\'ees de $ F$ \  \`a bords dans  $ F \cap \mathcal{U} $  modulo les bords  de $(n+1)$-cha\^ines compactes orient\'ees  et les n-cha\^ines compactes trac\'ees dans    $ F \cap \mathcal{U} $  (voir l'Annexe ).

\noindent Soit donc $\gamma $ \ une n-cha\^ine compacte orient\'ee de  $F $ \  \`a bord  dans  $ F \cap \mathcal{U} .$  Comme la monodromie  $T_{c}$  laisse  $ F \cap \mathcal{U} $ \ fixe, on aura  $ \partial \gamma = \partial T_{c} \gamma $ \ et \ $ (T_{c}-1) \gamma $  \ sera un $n$-cycle compact orient\'e de  $ F $.  Il est facile de voir que pour la somme d'un bord et d'une cha\^ine trac\'ee dans  $ F \cap \mathcal{U} $ , \ $ \gamma = \partial \Gamma + \Delta $ \ , on a  $ (T_{c}-1) \gamma = \partial  (T_{c}-1) \Gamma  $\ car \ $ (T_{c}-1) \Delta = 0 $ . Donc  $ var_{c} $ \ est bien d\'efinie et les relations annonc\'ees ci-dessus en d\'ecoulent imm\'ediatemment .

\noindent Pour calculer  $ var_{c} $ \  en cohomologie de de Rham, nous utiliserons le formalisme de [B.97]  qui s'appliquerait en fait au cas, plus g\'en\'eral , o\`u l'on aurait seulement l'hypoth\`ese  $ \bar{\mathcal{U}'} \cap S_{1} = \emptyset $ \ o\`u  $S_{1} $ \ d\'esigne la r\'eunion du lieu singulier de  $ {X}_{\mathbb{C}}$ \ et de l'ensemble des points lisses de  $ {X}_{\mathbb{C}}$ \ en lesquels la valeur propre  1  de la monodromie agissant sur la cohomologie r\'eduite de la fibre de Milnor de $f$  appara\^it ; on a donc  $ S_{1} \subset \widehat{S}  $ , in\'egalit\'e \'eventuellement stricte.

\noindent Introduisons alors  $$ \widetilde{var_{c}} :=  \Theta_{c} \circ var_{c} \ ,$$

\noindent o\`u $ \Theta_{c} $  est l'automorphisme de  $ H_{c}^{n}(F)_{=1}$ \ donn\'e par 
$$ \Theta_{c} : = \sum_{k=0}^{+ \infty}  \  \frac{(-1)^{k}}{k+1} (T_{c}-1)^{k} $$

\noindent On a alors   \ $  \widetilde{var_{c}} \circ can_{c} = \frac{i}{2 \pi} Log T_{c} $ \ sur  $ H_{c}^{n}(F)_{=1}$ \ . D\'efinissons l'application sesquilin\'eaire 
$$ h :   H_{c}^{n}(F,F\cap \mathcal{U})_{=1} \times H^{n}(F)_{=1}  \longrightarrow \mathbb{C} $$

\noindent en posant pour  $ (e,e') \in  H_{c}^{n}(F,F\cap \mathcal{U})_{=1} \times H^{n}(F)_{=1} $ :
$$  h(e,e') = \mathcal{I}( \widetilde{var_{c}} (e) , e') $$
o\`u $ \mathcal{I}  $ \ d\'esigne la dualit\'e de Poincar\'e hermitienne sur  $ F $  donn\'ee par  :
$$ \mathcal{I} :  H_{c}^{n}(F) \times H^{n}(F) \longrightarrow  \mathbb{C} \qquad  avec \qquad  \mathcal{I} (a,b) := \frac {1}{(2i \pi )^{n}} \int_{F}^{}  a \wedge \bar{b}   .$$

\noindent { \bf Remarque .}

\noindent Dans le cas o\`u $ \Lambda = f_{\mathbb{C}}^{-1}(0) \cap \partial B(0, \varepsilon '') $  , on  a une identification naturelle de   $H_{c}^{n}(F,F\cap \mathcal{U}) $  avec  $ H^{n}(F) $ \ et on retrouve la forme hermitienne canonique  : voir [B.85] dans le cas  $ X_{\mathbb{C}} $ lisse et  $f$ \ \`a singularit\'e isol\'ee , [B.90] pour le cas   $ X_{\mathbb{C}} $ lisse et  $f$ \ \`a singularit\'e  isol\'ee  pour la valeur propre 1  de la monodromie et  [B.M.02] pour le cas o\`u   $ X_{\mathbb{C}} $ \ et \ $f$ \  sont \`a singularit\'es  isol\'ees .

\bigskip

\section{\bf Le cas d'une valeur propre  $ \not= 1$ }

\noindent Pla{\c c}ons-nous dans la situation  pr\'ec\'edente
 en supposant de plus que  $(X_{\mathbb{C}}, 0)$  est le complexifi\'e d'un germe d'ensemble analytique r\'eel   $(X_{\mathbb{R}}, 0)  \subset ( \mathbb{R}^{N}, 0)$  admettant une singularit\'e isol\'ee en  $0$ .
\noindent Supposons de plus que le germe  $( f, 0) $ soit le complexifi\'e d'un germe analytique r\'eel  $(f_{\mathbb{R}}, 0) : (X_{\mathbb{R}}, 0) \longrightarrow  (\mathbb{R}, 0) $  admettant une singularit\'e isol\'ee en  $0$  sur  $(X_{\mathbb{R}}, 0)$ .

\noindent On notera   $ X_{\mathbb{R}} :=  X_{\mathbb{C}} \cap  \mathbb{R}^{N} $ \ et  $ f_{\mathbb{R}} $  la restriction de  $ f = f_{\mathbb{C}} $  \`a   $ X_{\mathbb{R}} $  \ , o\`u  $ X_{\mathbb{C}} $ d\'esigne un repr\'esentant de Milnor  du germe  $(X_{\mathbb{C}}, 0)$  comme pr\'ec\'edemment .
\noindent Nous supposerons dans tout ce qui suit que la vari\'et\'e analytique r\'eelle (lisse)   $ X_{\mathbb{R}} - \lbrace 0 \rbrace $  est orientable, et que nous avons fix\'e une orientation .

\noindent Nous fixerons \'egalement  $ 0 < \varepsilon' < \varepsilon'' < \varepsilon  $  \  avec  \  $ \varepsilon  $ \ fix\'e comme pr\'ec\'edemment  et \ $  \varepsilon -\varepsilon'  \ll \varepsilon $ \ . Nous supposerons que \  $ \varepsilon $ \   a \'et\'e choisi assez petit pour que  \ $ \Lambda = f_{\mathbb{R}}^{-1}(0) \cap \partial B(0, \varepsilon'') $ \  soit une sous-vari\'et\'e analytique r\'eelle lisse et compacte de  $ f_{\mathbb{C}}^{-1}(0) - \widehat{S }$ . Ceci est possible gr\^ace \`a nos hypoth\`eses de singularit\'es isol\'ees pour  $ X_{\mathbb{R}} $ \ et \ $ f_{\mathbb{R}}$ .

\noindent Pour chaque composante connexe  $ A_{\alpha} $ \ de  $ X_{\mathbb{R}} - f_{\mathbb{R}}^{-1}(0) $ (il n'y en a qu'un nombre fini d'apr\`es Lojasiewicz , car c'est la diff\'erence de deux ensembles semi-analytiques  compacts  , voir [Lj.65] ) on d\'efinit une $n-$cha\^ine compacte orient\'ee \   $ \Gamma(A_{\alpha}) $  \ de la fibre de Milnor  $ F $ , dont le bord est contenu dans  $ F \cap \mathcal{U} ,$ \ de la facon suivante  : 

\noindent Pour  $ A_{\alpha} \subset \lbrace f_{\mathbb{R}} > 0  \rbrace $ \ on pose 

$$ \Gamma(A_{\alpha}) :=   \Gamma(A_{\alpha})^{+} = A_{\alpha} \cap f_{\mathbb{R}}^{-1}(s_{0}) \cap \overline{B(0, \varepsilon'')} $$
o\`u  \  $ s_{0} \in D \cap \mathbb{R^{+*}} $ \ est le point base de \  $D^{*} = D_{\delta}^{*}  \ $  \ (rappelons que , par d\'efinition , la fibre de Milnor  $ F $  de  $ f_{\mathbb{C}} $ \ en  $0$  est \'egale \`a  \  $ f_{\mathbb{C}}^{-1}(s_{0}) $ )  . On oriente   \ $ f_{\mathbb{R}}^{-1}(s_{0}) $ \   comme le bord de l'ouvert  \ $ f_{\mathbb{R}}^{-1}( s < s_{0} ) $ \ de  \  $ X_{\mathbb{R}} - \lbrace 0 \rbrace $ . Donc  \  $ \Gamma(A_{\alpha}) $ \ d\'efinit une classe dans  \  $ H_{c}^{n}(F,F\cap \mathcal{U}) $ .

\noindent Pour  $ A_{\alpha} \subset \lbrace f_{\mathbb{R}} < 0  \rbrace $ \  on d\'efinit   \ $\Gamma(A_{\alpha}) $ \  \`a partir de la $n-$cha\^ine compacte ,\`a bord dans  $ \mathcal{U} $  :
$$ \Gamma(A_{\alpha})^{-} = A_{\alpha} \cap f_{\mathbb{R}}^{-1}(-s_{0}) \cap \overline{B(0, \varepsilon'')} $$
que l'on oriente  gr\^ace \`a l'orientation de  \ $ f_{\mathbb{R}}^{-1}(-s_{0}) $ \ comme bord de l'ouvert  \ $ f_{\mathbb{R}}^{-1}( s >- s_{0} ) $ \ de  \  $ X_{\mathbb{R}} - \lbrace 0 \rbrace $ , 
 en la suivant dans la trivialisation  $ \Psi $ \ de  $ f_{\mathbb{C}}$ \  le long du demi-cercle   \ $ \lbrace s_{0}e^{i\theta} , \theta \in [-\pi,0 ] \rbrace $ \ et en changeant l'orientation .

\noindent Dans ces conditions on notera  :
$$ \Gamma(A_{\alpha}) = - T^{\frac{1}{2}}. \Gamma(A_{\alpha})^{-}\ . $$ 

\noindent C'est \`a nouveau une  $n-$cha\^ine compacte orient\'ee contenue dans  $F $ et \`a bord dans  $ F \cap \mathcal{U} $ \ . Elle d\'efinit donc une classe dans  \  $ H_{c}^{n}(F,F\cap \mathcal{U}) $ .

\noindent Le changement d'orientation dans la d\'efinition de  \ $ \Gamma(A_{\alpha}) $ \ dans le cas n\'egatif est d\^u au fait que la rotation  de  $ \pi $  pour  $ \mathbb{R}^{-} $ \  l'envoie sur 
$  \mathbb{R}^{+} $ \ avec l'orientation "oppos\'ee" .

\noindent On \'etend alors l'application  $ \Gamma $ \  par  $ \mathbb{C}$-lin\'earit\'e en une application  $ \mathbb{C}$-lin\'eaire 
$$ \Gamma :  H^{0}(X_{\mathbb{R}}- f_{\mathbb{R}}^{-1}(0) ,\mathbb{C} ) \longrightarrow   H_{c}^{n}(F,F\cap \mathcal{U}) .$$

\noindent Si  \  $ A = \sum_{\alpha}^{ } a_{\alpha}.A_{\alpha} $ \ avec  \  $ A^{+} = \sum_{A_{\alpha} \subset \lbrace f_{\mathbb{R}}>0 \rbrace}^{ } a_{\alpha}.A_{\alpha} $ \ et  \ $ A^{-} = \sum_{A_{\alpha} \subset \lbrace f_{\mathbb{R}}<0 \rbrace}^{ } a_{\alpha}.A_{\alpha} $ \ on aura  :  $$  \Gamma(A) =  \Gamma(A)^{+}   - T^{\frac{1}{2}}. \Gamma(A)^{-} .$$ 

\noindent{\bf Th\'eor\`eme 1}

\bigskip

\noindent \textit{On se place sous les hypoth\`eses ci-dessus , \`a savoir que  $$  \widetilde{f_{\mathbb{R}}} :  (X_{\mathbb{R}}, 0) \longrightarrow  (\mathbb{R}, 0) $$  est un germe de fonction analytique r\'eelle \`a singularit\'e isol\'ee en $0$   sur un germe d'ensemble analytique r\'eel de dimension pure $n+1$  dans $ (\mathbb{R}^{N}, 0) $ , ayant une singularit\'e isol\'ee en  $0$ , v\'erifiant les hypoth\`eses  H a),H b),H c) et H d) du paragraphe  2 .
 On note par  $$f_{\mathbb{R}} : X_{\mathbb{R}}  \longrightarrow  ]- \delta, \delta[  $$  la trace sur le r\'eel d'un repr\'esentant de Milnor du complexifi\'e  $$  \widetilde{f_{\mathbb{C}}} :  (X_{\mathbb{C}}, 0) \longrightarrow  (\mathbb{C}, 0) $$ de  $ \widetilde{f_{\mathbb{R}}} $ . On d\'esignera par  $F$  la fibre de Milnor de   $\widetilde{f_{\mathbb{C}}}$ en \ $0$ } .

\noindent  \textit{Pour tout  $ u \in ]0,1[ $ \ l'application canonique 
$$ can :   H^{n}(F,F\cap \mathcal{U})_{e^{-2i\pi u}}  \longrightarrow  H^{n}(F)_{e^{-2i\pi u}} $$
est bijective , o\`u  \ $\mathcal{U} $  d\'esigne un voisinage ouvert convenable (voir plus haut) de \ $ \Lambda =  f_{\mathbb{R}}^{-1}(0) \cap \partial B(0, \varepsilon '') $ } \ .

\noindent \textit{Pour tout  $ A \in  H^{0}(X_{\mathbb{R}}- f_{\mathbb{R}}^{-1}(0) ,\mathbb{C} )  $ \ on aura : 
\begin{equation}
 (-2i\pi)^{n}.\langle \Gamma(A) ,\overline{can^{-1}(e) }\rangle  =  Res( \lambda = -u , \int _{A}^{}  (f/s_0)^{\lambda} \  \frac{df}{f} \wedge w_{k} ) 
\end{equation}
o\`u l'application \   $  \Gamma :   H^{0}(X_{\mathbb{R}}- f_{\mathbb{R}}^{-1}(0) ,\mathbb{C} ) \longrightarrow H^{n}(F,F\cap \mathcal{U}) $  est d\'efinie ci-dessus , o\`u  $ e \in H^{n}(F)_{e^{-2i\pi u}} $ \ , o\`u  $ \langle , \rangle $ \ d\'esigne l'application sesquilin\'eaire 
$$  \langle , \rangle  :  H_{c}^{n}(F,F \cap \mathcal{U})\times   H^{n}(F,F \cap \mathcal{U}) \rightarrow  \mathbb{C}  $$ 
d\'eduite de la dualit\'e de Poincar\'e sur  $F$ (voir l'Annexe ) , et o\`u  $ w_{1},\ldots,w_{k} $ \ sont des $n-$formes semi-m\'eromorphes \`a p\^oles dans
$ f_{\mathbb{C}}^{-1}(0) $ \ sur $ X_{\mathbb{C}} $ \ v\'erifiant  les conditions   suivantes  : 
\begin{equation}
 \begin{cases}  
$$ dw_{j} = u\frac{df}{f} \wedge w_{j} + \frac{df}{f} \wedge w_{j-1} \ \  \  \  \forall j\in [1,k] ,avec  \  \ \  w_{0} = 0 $$ \\
$$ \  w_{k}\vert_{F} = e  \  \  \ dans  \    \  \  H^{n}(F)   \qquad  \qquad  \qquad  \qquad  \qquad\  \Box $$ 
\end{cases}
\end{equation}}


\noindent { \bf Remarques :}
\begin{itemize}
\item[1)]  Rappelons que  $ F := f_{\mathbb{C}}^{-1}(s_{0})  $ \ o\`u \ $ s_{0}\in D \cap \mathbb{R}^{+*}  $ \ est le point base ; donc  $ w_{k} $ \  est une forme  $ \mathcal{C}^{\infty}  $ \ et $d-$ferm\'ee de degr\'e n  sur  $F$ .
\item[2)] Dans la formule du th\'eor\`eme on a omis une fonction de troncature  $ \rho \in \mathcal{C}_{c}^{\infty}(X_{\mathbb{R}}) $ \  valant identiquement  1  pr\`es de  0  qui est n\'ecessaire pour donner un sens \`a l'int\'egrale  $ \int _{A}^{}  \rho . (f/s_0)^{\lambda} \  \frac{df}{f} \wedge w_{k} $ \ ; en effet elle n'est pas utile pour discuter les p\^oles de  $ \frac{1}{\Gamma(\lambda)} \int _{A}^{}  (f/s_0)^{\lambda} \  \frac{df}{f} \wedge w_{k} ) $ \ qui sont concentr\'es \`a l'origine  et dont les parties polaires sont des \ $(n+1)-$courants  \`a support  $ \lbrace 0 \rbrace $ \ vu nos hypoth\`eses  \qquad \  \qquad \qquad \qquad $ \Box $ 

\end{itemize}

\bigskip

\noindent{\bf Corollaire : }

\noindent \textit{Sous les hypoth\`eses du Th\'eor\`eme 1  l'ordre des p\^oles du prolongement m\'eromorphe de  $ \int_{A}^{} f^{\lambda} \  \Box $ \  en  \ $ -u-\nu , \nu \in \mathbb{N} , \nu \gg 1  ,$ \ est exactement l'ordre de nilpotence de  $ T - e^{-2i\pi u}   $ \ agissant sur la composante de  $ \Gamma(A) $  dans  $  H^{n}(F)_{e^{-2i\pi u}}  \qquad  . $ }

\noindent \textit{En particulier , la composante de   $ \Gamma(A) $  dans  $  H^{n}(F)_{e^{-2i\pi u}} $ \ est nulle si et seulement si le prolongement m\'eromorphe de  $ \int_{A}^{} f^{\lambda}  \ \Box $ \ n'a jamais de p\^ole en  $ -u-\nu , \forall \nu \in \mathbb{N}  \qquad  \qquad  \qquad  \qquad  \Box $}

\bigskip

\bigskip

\noindent Preuve du th\'eor\`eme 1 :

\noindent Avant de passer \`a la preuve proprement dite, il est bon de pr\'eciser la normalisation que nous utilisons ici pour repr\'esenter un \'el\'ement de   $  H^{n}(F)_{e^{-2i\pi u}} .$ En effet, cette normalisation  est celle d\'ej\`a utilis\'ee dans  [B.M.02] et diff\`ere de celle de  [B.91] . Comme les formules analogues \`a celle du th\'eor\`eme 1  ci-dessus sont erron\'ees, commen{\c c}ons par discuter ces deux normalisations et corriger les formules de  [B.M.02] .

\bigskip

\noindent { \bf Remarque/Erratum.}

\bigskip

\noindent  Soit  $ u \in [0,1[ $  et  soit  $ e\in H^{n}(F)_{e^{-2i\pi u}} .$  
Pour "repr\'esenter" la classe de cohomologie \ $  e$ \  par une section m\'eromorphe (uniforme) du fibr\'e de Gauss-Manin de  f   on utilise , dans  [B.91]  par exemple, la normalisation suivante :

\noindent   Soit 
 $$ \mathcal{N} := \frac{i}{2\pi}Log(e^{-2i\pi u}.T\vert_{H^{n}(F)_{e^{-2i\pi u}}  }) $$
 ou  $-2i\pi  \mathcal{N} $ est le logarithme nilpotent de l'endomorphisme unipotent 
$$e^{-2i\pi u}.T\vert_{H^{n}(F)_{e^{-2i\pi u}}  }. $$

\noindent  On pose  
$$\varepsilon := exp[(u + \mathcal{N})Log s](e)  = s^{u}.\sum_{0}^{\infty} \frac{(Log s)^{j}}{j!} \mathcal{N}^{j}(e) . $$

\noindent Bien s\^ur, la monodromie qui agit par  $ Log s  \longrightarrow   Log s + 2i\pi $ \   et  par  \ 
 $ e  \longrightarrow   T(e) ,$ laisse  $ \varepsilon$  invariante ; et , puisque la fibre de Milnor est d\'efinie par  $  F := f^{-1}(s_0) $  o\`u  $ s_0 \in  D\bigcap \mathbb{R}^{+*}  $ est le point base choisi ,
on aura  \  
$$ \varepsilon \vert_{F} =  s_{0}^{u}.\sum_{0}^{\infty} \frac{(Log s_{0})^{j}}{j!} \mathcal{N}^{j}(e) $$
  dans $H^{n}(F)_{e^{-2i\pi u}}  $.
Pour retrouver \ $ e$ \   \`a partir de  $\varepsilon $  on utilise le morphisme  :

$$ r^{n}(k)  :   h^{n}(k) _{0}    \longrightarrow   Ker \  \mathcal{N}^{k} \subset H^{n}(F)_{e^{-2i\pi u}} $$
o\`u  $ h^{n}(k)$  est le i-\`eme faisceau de cohomologie du complexe  
$ (\Omega ^{\cdot}(k),\delta_{u}) $  (voir  [B.91] ).

\noindent  Si l'on pose   $\varepsilon _{j} :=  \mathcal{N}^{k-j}(\varepsilon) $ \  pour \ $  j\in [1,k]$ , alors
$ \breve{ \varepsilon} :=  \lbrace  \varepsilon_{k},..., \varepsilon_{1} \rbrace $ \  est une section de \  $ h^{n}(k)$  (car  $\delta_{u}(\breve{ \varepsilon}) = 0 $)  et on a
  $$ r^{n}(k)( \breve{ \varepsilon} ) =  e  . $$ 

\noindent Nous allons d\'ecrire maintenant une autre normalisation qui permet de simplifier les formules  ; posons   
$$ \tilde{\varepsilon}:= exp[(u + \mathcal{N})Log( s/s_{0})](e) . $$ 

\noindent La section m\'eromorphe du fibr\'e de Gauss-Manin donn\'ee par  $ \tilde{\varepsilon}$  est reli\'ee \`a   \ $\varepsilon $ \  par  l'\'egalit\'e  
$$exp[(u + \mathcal{N})Log (s_{0})](\tilde{\varepsilon}) = \varepsilon $$

\noindent et elle v\'erifie simplement que   $ \tilde{\varepsilon} \vert_{F} =  e $ .

\noindent Le prix \`a payer pour obtenir des formules plus simples gr\^ace \`a cette seconde "normalisation"  est de remplacer  $ f^{\lambda} $  \  par  $ (f/s_0)^{\lambda} $ \  dans les int\'egrales consid\'er\'ees.

\noindent On prendra garde au fait que dans  [B.M.02]  les formules  6.2a, 6.2b et  6.9 sont erron\'ees : on doit   remplacer  $ f^{\lambda} $  \  par  $ (f/s_0)^{\lambda} $ \  puisqu'on  y a adopt\'e la seconde normalisation d\'ecrite ci-dessus  ( donc  $ \varepsilon_{k} \vert_{F} = e $ \ ) , le coefficient
 $ s_{0}^{r} $  \ de la formule  6.2a  devant \^etre supprim\'e .

\noindent En fait l'erreur de calcul dans  [B.M.02]  se trouve apres la formule  6.4 .Il faut rectifier de la fa{\c c}on suivante  :

 \begin{eqnarray*}
 \int_{A}( f/s_{0})^{\lambda} w_{k}\wedge \frac{df}{f} \qquad \qquad \qquad  =  \qquad \qquad \qquad &  \qquad \qquad \qquad  & \\
\int _{0}^{s_{0}}(s/s_{0})^{\lambda} \ a(s) \frac{ds}{s}   -  \int _{0}^{s_{0}} (s/s_{0})^{\lambda}\ H(s) \frac{ds}{s}  &  (mod. \  ent. funct.)  & \\  =  \int _{0}^{s_{0}} (s/s_{0})^{\lambda + r}\  [P - Q](Log( s/s_{0})) \frac{ds}{s} &  (mod. ent. funct.) & 
\end{eqnarray*}

\noindent In consequence, \ the residue at  \  $ \lambda = -r $ \  of \  $ \int_{A}( f/s_{0})^{\lambda} w_{k}\wedge \frac{df}{f} $ \ is equal   to  \  P(0) - Q(0)  \ $ \ \ \ \ \ \blacksquare $

\bigskip

\noindent Revenons \`a la preuve du th\'eor\`eme 1 . La bijectivit\'e de l'application canonique 
$$ can :   H^{n}(F,F\cap \mathcal{U})_{e^{-2i\pi u}}  \longrightarrow  H^{n}(F)_{e^{-2i\pi u}} $$
a \'et\'e \'etablie au paragraphe  3  .
\noindent Montrons d\'ej\`a que le second membre de la formule du th\'eor\`eme ne depend pas du choix de   $ w_{1},\ldots,w_{k} $  v\'erifiant les conditions (2)  
de l'\'enonc\'e . On sait , d'apr\`es [B.91] par exemple , que si  $ w' = (w'_{1},\ldots,w'_{l}) $ repr\'esente aussi $  e$  \ on peut trouver des $n-$formes semi -m\'eromorphes $ \alpha$ \  et  $  \beta$ \ \`a p\^oles dans  $ f^{-1}(0) $ \  telles que l'on ait  
$$  w_{k} - w'_{l} = d\alpha + df\wedge \beta $$
 \ sur  $ X_{\mathbb{C}} $ .

\noindent Il s'agit alors de voir que le prolongement m\'eromorphe de
$$ \int_{A}^{} f^{\lambda}\rho \frac{df}{f} \wedge d\alpha  $$ n'a pas de p\^ole congrus a \  $-u $ \ modulo $ \mathbb{Z} $ . Or la formule de Stokes donne pour $ \Re(\lambda) \gg 1 $  \ puisque  \ $ \partial A \cap Supp(\rho) $ \   est contenu dans  $ f^{-1}(0) $ :
$$  \int_{A}^{} f^{\lambda}\rho \frac{df}{f} \wedge d\alpha =   - \int_{A}^{}d ( f^{\lambda}\rho \frac{df}{f} \wedge \alpha) +  \int_{A}^{} f^{\lambda}d\rho \wedge \frac{df}{f} \wedge \alpha $$ 
et donc  
$$ \frac{1}{\Gamma(\lambda)} \int_{A}^{} f^{\lambda}\rho \frac{df}{f} \wedge d\alpha $$
 n'aura pas de p\^ole du tout puisque  $ d\rho \equiv 0 $ \  pr\`es de  $0$  .

\noindent Nous allons maintenant choisir des repr\'esentants  $ \check{w}_{1}, \ldots, \check{w}_{k} $ \ v\'erifiant la condition suppl\'ementaire  :
$$  Supp(\check{w_{j}}) \cap \mathcal{U} = \emptyset   \  \  \forall j \in [1,k]  .$$

\noindent Pour cela nous allons utiliser l'annulation du groupe d'hypercohomologie 
$$  \mathbb{H}^{n}(\mathcal{U}',\mathcal{E}^{\cdot}(k),\delta^{\cdot}_{u} )  =  0  $$ 
qui r\'esulte du fait que  $ \mathcal{U}' \cap \widehat{S} = \emptyset $ \ et que les faisceaux de cohomologie  $ h^{i}(k) $ \ du complexe  $ (\mathcal{E}^{\cdot}(k),\delta^{\cdot}_{u} ) $ \ sont support\'es par  $ \widehat{S }$  (car on suppose  $ e^{-2i\pi u}  \not= 1 $ , voir  [B.91] prop.1 p.427  ) .

\noindent On peut donc trouver  $ \xi \in H^{0}(\mathcal{U}',\mathcal{E}^{n-1}(k))  $ \  v\'erifiant  :
$$  \delta_{u} \xi  =  w \vert_{\mathcal{U}'}  $$
si   $ w \in H^{0}(X,\mathcal{E}^{n}(k)) $ \  v\'erifie les conditions  (2) de l'\'enonc\'e du th\'eor\`eme .
Choisissons alors une fonction  $  \sigma \in \mathcal{C}^{\infty}(X) $ \  v\'erifiant  \ $ \sigma \equiv 1  $ \ sur $  \mathcal{U} $ \  et  $ \ Supp(\sigma) \subset \mathcal{U}'  $ . Posons alors  :
$$  \check{w} :=  w -  \delta_{u}(\sigma.\xi)  .$$
Alors  $ \check{w} $ \ v\'erifie encore les conditions  (2)  et satisfait la condition de support d\'esir\'ee .
On remarquera qu'alors la restriction de  $ \check{w} $   \`a  $ X_{\mathbb{R}} \cap B(0,\varepsilon'') $ \ a un support  $ f_{\mathbb{R}}-$propre . 

\noindent Si maintenant la fonction  $ \rho \in \mathcal{C}^{\infty}(X_{\mathbb{R}})$  vaut identiquement  1  sur  $ B(0,\varepsilon'')  $\  et a son support contenu dans    $ (X_{\mathbb{R}} \cap B(0,\varepsilon'') )\cup \mathcal{U} $ \ , la fonction m\'eromorphe 
 $$\frac{1}{i\pi} \int_{A}^{} f^{\lambda} \rho \frac{df}{f}\wedge \check{w}_{k}  $$
sera la transform\'ee de Mellin , \`a une fonction enti\`ere pr\`es , de la fonction 
$$  s \longrightarrow  \int_{(f_{\mathbb{R}}= s)\cap A }\check{w}_{k}   .$$
On conclut alors la preuve du th\'eor\`eme 1 comme dans  le paragraphe 6 , Th.6.1 a)   de  [B.M.02] \  (modulo l'erratum ci-dessus )   \  \   \qquad \qquad \qquad $\blacksquare$ 

\bigskip

\noindent Preuve du Corollaire :

\noindent Un calcul \'el\'ementaire montre, qu'avec les notations introduites ci-dessus on a
$$  ( s \frac{\partial}{\partial s} -  u )^{k-1} ( \int^{}_{(f_{\mathbb{R}}=s) \cap A} \check{w_{k}} \ )  =  \int^{}_{(f_{\mathbb{R}}=s) \cap A} \check{w_{1}}  $$
et donc l'ordre de nilpotence de  $ \Gamma(A) $  dans \  $  H^{n}(F)_{e^{-2i\pi u}} $ \  minore l'ordre des p\^oles  en  $ -u -\nu , \nu \in \mathbb{Z} $ \ du prolongement m\'eromorphe de  $ \int_{A}^{} f^{\lambda} \Box . $ \ En effet, comme l'application canonique :
$$ can_{c}:   H_{c}^{n}(F)_{e^{-2i\pi u}} \longrightarrow  H_{c}^{n}(F,F\cap \mathcal{U})_{e^{-2i\pi u}}  $$ est bijective, puisque la monodromie agit comme l'identit\'e sur  $  H_{c}^{*}(F\cap \mathcal{U}) $  , la dualit\'e de Poincar\'e (hermitienne) entre  $  H_{c}^{n}(F)_{e^{-2i\pi u}} $ \ et \ $ H^{n}(F)_{e^{-2i\pi u}} $ \ montre que l'accouplement sesquilin\'eaire  :
 $$ H_{c}^{n}(F,F \cap \mathcal{U}))_{e^{-2i\pi u}}\times   H^{n}(F)_{e^{-2i\pi u}} \rightarrow \mathbb{C}  $$  est une dualit\'e hermitienne  pour laquelle la monodromie est auto-adjointe, d'o\`u notre assertion .

\noindent Pour voir que l'ordre des p\^oles ne peut d\'epasser l'ordre de nilpotence de  $ \Gamma(A) $ \ dans  $  H_{c}^{n}(F,F\cap \mathcal{U})_{e^{-2i\pi u}}  $ \ , on utilise le fait que les parties polaires des p\^oles congrus \`a \ $ -u $ modulo $\mathbb{Z} $  \ sont des distributions \`a support l'origine. Elles sont donc d'ordre fini , et pour calculer l'une d'entre elles, on peut remplacer toute forme-test par son d\'eveloppement de Taylor en  $0$ \`a un ordre assez \'elev\'e  (fixe, ne d\'ependant que de la partie polaire consid\'er\'ee ). On se ram\`ene alors au cas de la restriction \`a  $ X_{\mathbb{R}} $  d'une  $(n+1)-$forme m\'eromorphe sur  $X_{\mathbb{C}} $  que l'on d\'ecompose alors dans le syst\`eme de Gauss-Manin localis\'e en  $ f_{\mathbb{C}} . $ \  On conclut alors facilement (pour ce type de raisonnement voir  [B.85] ou  [B.M.00]  pour plus de d\'etails )   \qquad   \qquad       $\blacksquare$ 

\bigskip

\newpage

\section{\bf Le cas de la valeur propre  1.}

\noindent Traitons maintenant le cas de la valeur propre 1 

\bigskip

\bigskip

\noindent{\bf Th\'eor\`eme 1 bis} :

\noindent \textit{ Dans la situation du th\'eor\`eme 1 , pour tout   $ A \in  H^{0}(X_{\mathbb{R}}- f_{\mathbb{R}}^{-1}(0) ,\mathbb{C} )  $ \ et tout  $ e \in H^{n}(F)_{=1} $ on a 
\begin{equation}
 (-2i\pi)^{n} h(\Gamma(A),\overline{e}) = Res(\lambda = 0,\frac{1}{\Gamma(\lambda)}\int_{A}^{} (f/s_{0})^{\lambda} \frac{df}{f} \wedge w_{k} )  
\end{equation}
o\`u   $ w_{1},\ldots,w_{k} $ \ v\'erifient les conditions  (2) du th\'eor\`eme 1  avec  $ u=0 $  , et o\`u  h  est la forme sesquilin\'eaire 
$$ h :   H_{c}^{n}(F,F\cap \mathcal{U})_{=1} \times H_{c}^{n}(F)_{=1}  \longrightarrow \mathbb{C} $$
d\'efinie au paragraphe 3  $  \qquad \qquad \qquad \Box  $ } 

\bigskip

\noindent{\bf Corollaire} :

\noindent \textit{Dans la situation du th\'eor\`eme 1 bis  l'ordre des p\^oles du prolongement m\'eromor\-phe de la distribution  $ \frac{1}{\Gamma(\lambda)}\int_{A}^{} f^{\lambda} \ \Box  $  aux entiers n\'egatifs assez grands (en valeur absolue  ) est \'egal \`a l'ordre de nilpotence de la monodromie agissant sur \  $ var_{c}(\Gamma(A)_{=1}) \in  H_{c}^{n}(F)_{=1} \  \qquad  \  \qquad  \  \Box $  }

\bigskip

\noindent Preuve du Th\'eor\`eme 1 bis:

\noindent Comme on sait d\'ej\`a  que la dualit\'e de Poincar\'e hermitienne  $ \mathcal{I} $  \ ainsi que l'accouplement sesquilin\'eaire  
$$ \langle , \rangle  :  H_{c}^{n}(F,F \cap \mathcal{U})\times   H^{n}(F,F \cap \mathcal{U}) \longrightarrow  \mathbb{C}  $$
d\'ecrit  dans l'Annexe , sont non d\'eg\'en\'er\'ees , d\'efinissons l'application de variation  :
$$ var  :   H^{n}(F)_{=1} \longrightarrow   H^{n}(F,F \cap \mathcal{U})_{=1}  $$
en posant  $  \langle \varepsilon, var(e) \rangle : = \mathcal{I}( var_{c}(\varepsilon), e)  $ .

\noindent Posons alors  $  \widetilde{var} : = \Theta \circ var  , $ \ o\`u \ $ \Theta $\ est l'analogue sur \ $ H^{n}(F)_{=1} $ \ de \ $ \Theta_{c} $ \ introduit au paragraphe  3 .

\noindent Nous allons v\'erifier que cette d\'efinition est bien compatible avec la description de  $ \widetilde{var} $ \ qui est donn\'ee  (en cohomologie de de Rham) dans  [B.97] . Cette v\'erification est particuli\`erement bien venue dans ce texte puisqu'elle repose fondamentalement sur le "Residu de J.Leray " !

\noindent Pour cela, fixons  $[\gamma] \in  H_{c}^{n}(F,F \cap \mathcal{U})_{=1} $ \ o\`u  $\gamma$ \ est une $n-$cha\^ine compacte orient\'ee de  $F$  \`a bord dans  $ F\cap \mathcal{U} $ \ et soit  $e$ \ un vecteur propre de la monodromie dans  $   H^{n}(F)_{=1} $ . En utilisant l'invariance de ces variations  par la monodromie, il est facile de voir que l'on peut se ramener \`a ce cas pour faire notre v\'erification  (\`a l'aide d'une r\'ecurrence sur l'ordre de nilpotence ). Soit donc  $w $ \ une $n-$forme semi-m\'eromorphe sur  $X_{\mathbb{C}} $ \ , \`a p\^oles dans  $ f^{-1}(0) $ \ , v\'erifiant  $ dw = 0 $ \ et \ $ w \vert_{F} = e  .$ 
Soit  $\sigma $  une fonction de  $ \mathcal{C}^{\infty}_{c}(X_{\mathbb{C}}) $ \ identiquement \'egale \`a  1  sur un voisinage de  $ X_{\mathbb{C}} - \mathcal{U}' $ \ et identiquement nulle sur  $\mathcal{U} $ . Alors la classe  $[\gamma] \in  H_{c}^{n}(F,F \cap \mathcal{U})_{=1} $ \ est repr\'esent\'ee par le courant  $ \sigma.\gamma  . $
\noindent Appliquons le r\'esultat de J. Leray  [L.59]  \`a  \ $ w\vert_{\mathcal{U'}} $ \ . On peut alors trouver des formes \ $ \mathcal{C}^{\infty}(\mathcal{U'}) , \ \eta_{1} $ \ de degr\'e \ $ n-1 $ \ ,\ $ \omega_{1} $ \ de degr\'e  $n$ \ et d-ferm\'ees  ainsi qu'une forme semi-m\'eromorphe \`a p\^oles dans \ $ f^{-1}(0) , \ \alpha_{1} $ \  de degr\'e \ $n-1$ \ v\'erifiant  :
\begin{equation}   w =  \frac{df}{f} \wedge \eta_{1}  + \omega_{1}  + d\alpha_{1}   \ .
\end{equation} 
Posons alors \ $ \eta = -d\sigma \wedge \eta_{1} ,\ \omega = d\sigma \wedge \omega_{1} $ \ et \ $ \alpha = \sigma.d\alpha_{1} $ \ . On obtient alors 
\begin{equation}  d\sigma \wedge w =  \frac{df}{f} \wedge \eta  + \omega  + d\alpha   \ .
\end{equation} 
o\`u \ $\eta$ \ est une \ $n$-forme dans \ $ \mathcal{C}^{\infty}(\mathcal{U'}) $ \ \`a support dans \ $ \mathcal{U'} - \mathcal{U} $\ o\`u \ $\omega $ \ est une \ $(n+1)$-forme dans \ $ \mathcal{C}^{\infty}(\mathcal{U'}) $ \ \`a support dans \ $ \mathcal{U'} - \mathcal{U} $\ avec \ $ d\eta = 0 $ \ et
\ $d\omega = 0 $ \ et o\`u \ $\alpha $ \ est semi-m\'eromorphe de degr\'e \ $n$ \ \`a p\^oles dans \ $f^{-1}(0)$ \ et \`a support dans \  $ \mathcal{U'} - \mathcal{U} $ .
D'apr\`es  [B.97] la classe \ $[\eta] \in  H^{n}(F,F \cap \mathcal{U}) $ \ repr\'esente \ $var_{1}(e) $ \ c'est \`a dire la variation telle qu'elle est d\'efinie dans [B.97]  (et m\^eme \ $\widetilde{var_{1}(e)} $ \ puisque  $ T(e) = e $ ).

\noindent Consid\'erons maintenant les fonctions 
$ s \longrightarrow \int_{\gamma_{s}}^{} w  $  \  et \   $ s \longrightarrow \int_{\gamma_{s}}^{}\sigma .w  $ \ 
o\`u  $ \gamma_{s} $  est la famille horizontale de $n-$cha\^ines compactes \`a bord dans  $ \mathcal{U} $ \ que l'on d\'eduit de  $\gamma$  \`a l'aide de notre trivialisation  $ \Psi $  .

\noindent Ces deux fonctions diff\`erent d'une fonction semi-m\'eromorphe sur $ D $  \`a p\^ole en  $0$ \, puisque  $ \gamma - \sigma.\gamma $ \ est \`a support dans  $\mathcal{U}'$ \ et puisque la trivialisation  $ \Psi \vert{\mathcal{U}'} $ \ se prolonge de fa{\c c}on  $ \mathcal{C}^{\infty}$ \ \`a \ $ D $ .

\noindent Comme  $\sigma \equiv 0 $ \  pr\`es du bord de  $\gamma_{s} $ \ qui est dans  $ \mathcal{U} $ \ on aura  :
$$  \Omega : = d (\int_{\gamma_{s}}^{} \sigma .w ) = \int_{\gamma_{s}}^{} d\sigma \wedge w  $$
au sens de l'image directe des formes \`a support propre par une submersion  ; $ \Omega $ \ est donc semi-m\'eromorphe sur  $D$\ de degr\'e  1  \` a p\^ole en  $0$\   .

\noindent Maintenant la formule (5)  donne  :
$$  \Omega = \frac{ds}{s} \int_{\gamma_{s}}^{} \eta  \  +  \int_{\gamma_{s}}^{} \omega  +  d(\int_{\gamma_{s}}^{} \alpha )  $$
car on a  $ \int_{\gamma_{s}}^{} d\alpha   = d(\int_{\gamma_{s}}^{} \alpha )  $ \ d'apres Stokes  (rappelons que  $ Supp(\alpha)\cap \mathcal{U} = \emptyset $ \ donc on \'evite le bord de  $\gamma_{s} $ ) .

\noindent La fonction  $ s \rightarrow  \int_{\gamma_{s}}^{} \eta $ \ est dans \ $ \mathcal{C}^{\infty}(D) $ \ car  $ Supp(\eta) \cap \mathcal{U} = \emptyset $ \ et \ $ \eta $ \ est \ $ \mathcal{C}^{\infty} $ .

\noindent De m\^eme, la  1-forme sur $D$ \ $ \int_{\gamma_{s}} \omega $ \ est \ $ \mathcal{C}^{\infty} $ \ et d-ferm\'ee  car  $ Supp(\omega ) \cap \mathcal{U} = \emptyset $ \ et $ d\omega = 0 $ .

\noindent On peut donc appliquer le lemme suivant , que Jean Leray ne saurait renier 

\bigskip

\noindent {\bf Lemme }

\noindent \textit{Soit \ $\varphi \in \mathcal{C}^{\infty}(D) $ \ et \ $\omega \in \mathcal{C}^{\infty}(D)^{1} $ \ telles que  $ \Omega := \frac{dz}{z} \varphi  +  \omega  $ \ soit  d-ferm\'ee sur  $D^{*}$ . }

\noindent \textit{Alors pour tout chemin ferm\'e  $ C $ \  de $ D^{*}$ \ d'indice  1  par rapport \`a l'origine , on a  $$  \int_{C}^{} \Omega  =  2i\pi .\varphi(0)  \  \ . $$ }

\noindent  Preuve :

\noindent Comme le cas o\`u  $ \Omega = \varphi(0).\frac{dz}{z} $ \ est clair, il s'agit de voir que pour  $ \psi \in \mathcal{C}^{\infty}(D) $ \ et  \ $\omega \in \mathcal{C}^{\infty}(D)^{1} $ \ la forme \   $ \Omega_{1} = \frac{\bar{z}}{z} \psi dz  + \omega  $  \ est d-ferm\'ee  si et seulement elle est d-exacte sur  $ D^{*}$ . Mais la  d-fermeture de  $\Omega_{1}$ \ implique  $ \frac{\partial}{\partial \bar{z}}( \frac{\bar{z}}{z} \psi) \in  \mathcal{C}^{\infty}(D) $ \ et comme  $ \frac{\partial}{\partial \bar{z}} : \mathcal{C}^{\infty}(D)         \longrightarrow \mathcal{C}^{\infty}(D) $ \ est surjective, il existe  $ \eta \in  \mathcal{C}^{\infty}(D) $ \        telle que  $ \frac{\partial}{\partial \bar{z}}(\eta -  \frac{\bar{z}}{z} \psi ) = 0 $ \ sur  $ D^{*}$ .

\noindent Comme la fonction  $ \eta -  \frac{\bar{z}}{z} \psi  $ \ est holomorphe sur   $ D^{*}$  \ et born\'ee en  $0$ , elle est holomorphe sur  $D$ . On en d\'eduit que  $  \frac{\bar{z}}{z} \psi  $ \  se prolonge de fa{\c c}on  $  \mathcal{C}^{\infty} $ \ \`a  $D$ \ et donc aussi  $\Omega_{1} $ . Alors $\Omega_{1} $ \ est d-exacte  \   \   $ \qquad \qquad  \qquad \qquad  \blacksquare$ 

\bigskip

\noindent {\bf Remarque}

\noindent Si dans le lemme pr\'ec\'edent on remplace  $\Omega $ \ par  $ \Omega + dg $ \ ou  $ g \in  \mathcal{C}^{\infty}(D^{*}) $ , la conclusion est la m\^eme . En particulier ce sera le cas quand   $ g = \frac{\gamma}{z^{M}} $ \ avec  $ \gamma \in \mathcal{C}^{\infty}(D) $ , c'est a dire pour  g  semi-m\'eromorphe \`a p\^ole en  $0 \qquad \qquad   \qquad   \Box $

\bigskip

\noindent On obtient donc en appliquant ce lemme, que pour tout chemin ferm\'e  $  C $ de  $ D^{*}$ \ d'indice 1  par rapport \`a l'origine on aura  :
$$  \int_{C}^{} \Omega  =  2i\pi \int_{\gamma_{0}} \eta  $$
On prendra garde que  $ \gamma_{0} $  n'est bien d\'efini que dans  $\mathcal{U}' $  \ ; comme le support de  $\eta $ est contenu dans  $\mathcal{U}' $ , notre formule a  bien un sens .

\noindent La fonction  $ s \rightarrow \int_{\gamma{s}}^{} \eta $ \ est constante puisque  $\eta $ \ est d-ferm\'ee \`a support disjoint de  $\mathcal{U} $ \ . On a alors  , puisque  $\eta$ repr\'esente $var_{1}(e) $
telle qu'elle est definie dans  [B.97]
$$  \int_{\gamma_{0}} \eta  = \int_{\gamma} \eta  = \int_{\gamma} var_{1}(e)  $$
Mais on a \'egalement  :
$$  \int_{C}^{} \Omega  = 2i\pi \int_{var_{c}(\gamma)}^{} e  $$
puisque  $ var_{c}(\gamma) = T\gamma - \gamma $ . On en conclut que  :
$$ \mathcal{I} ( var_{c}(\gamma) ,\bar{e} ) = \langle \gamma , \overline{var_{1}(e)} \rangle  $$
ce qui montre que la variation $ var_{1} $ \ d\'efinie dans  [B.97] est bien "l'adjoint" de  $ var_{c} $ \ et coincide  avec la variation d\'efinie plus haut . Ceci justifie donc l'usage du  calcul de la variation en cohomologie de de Rham  de [B.97] .

\bigskip

\noindent Passons maintenant \`a la preuve proprement dite du th\'eor\`eme 1 bis .

\noindent Consid\'erons $ e \in H^{n}(F)_{=1} $ \ et \ $ w_{1},\ldots,w_{k} $ \ repr\'esentant  \ $e$ \ c'est \`a dire v\'erifiant les conditions (2) de l'\'enonc\'e du th\'eor\`eme 1  avec  $u = 0 $ \ . Comme la monodromie est l'identit\'e sur  $ \mathcal{U}' $  \ on peut  trouver \ $ \xi \in   H^{0}(\mathcal{U}',\mathcal{E}^{n-1}(k)) $ \ v\'erifiant sur   $ \mathcal{U}'$   \  :
$$  \delta_{0} \xi  =  \mathcal{N}_{k} w   \ \   $$
o\`u  $  \mathcal{N}_{k} $ d\'esigne l'endomorphisme du complexe 
 $ (\mathcal{E}^{\cdot}(k),\delta^{\cdot}_{0}) $ qui induit  $ \frac{i}{2\pi}Log T $  sur  \ $ H^{n}(F)_{=1} $ \ (voir [B.91] ou [B.97] ). Choisissons alors une fonction \ $  \rho$ \ dans \ $ \mathcal{C}^{\infty}(X) $ \  v\'erifiant  \ $ \rho \equiv 1  $ \ sur $  \mathcal{U}'' $ \  , un voisinage ouvert de \ $ \overline{\mathcal{U}} $ \ dans  $ \mathcal{U}' $  \  ,
 et  $ \ Supp(\rho) \subset \mathcal{U}'  $ \ .
Posons  (comparer avec [B.97] p.15 ) :
$$  v :=  \mathcal{N}_{k} w -  \delta_{0}(\rho \xi)   \  \  .$$
On a alors  \ $ \delta_{0} v = 0  \ , \  v\vert_{X -\overline{\mathcal{U''}}} =  \mathcal{N}_{k} w   $ \ et \ $ Supp(v)\cap \overline{\mathcal{U''}} = \emptyset $ .

\noindent Ceci conduit explicitement aux relations suivantes :
$$ w_{j-1} = v_{j} + d(\rho \xi_{j})  - \frac{df}{f} \wedge \rho \xi_{j-1} \  \  \  \    \forall  j \in [1,k]  $$
avec les conventions  $w_{0} = 0 \ , v_{0} = 0 \  , \xi_{0} = 0 $ \ .

\noindent Choisissons maintenant  une fonction  $  \sigma \in \mathcal{C}^{\infty}(X) $ \  v\'erifiant  \ $ \sigma \equiv 0  $ \ au voisinage de \ $ \overline{ \mathcal{U}}$ \ et \ $ \sigma \equiv 1  $ \ au voisinage de $ X -  \mathcal{U}''$ \ . La relation  $ \delta_{0}(w) = 0 $ \  qui est cons\'equence de la condition (2) pour  $u=0$ \ donne $$ dw_{k} = \frac{df}{f} \wedge w_{k-1} . $$ 
Posons  $  W := d\sigma \wedge ( w_{k} + \frac{df}{f} \wedge \rho \xi_{k}) $ \ . C'est une forme semi-m\'eromorphe d-ferm\'ee et \`a support  f-propre dans  $ \mathcal{U}' - \overline{\mathcal{U}} $ \ . Soit  $ Y' := f^{-1}(0) \cap ( \mathcal{U}' -  \overline{\mathcal{U}} ). $ \ Comme  $ \widehat{S}\cap  \mathcal{U}' = \emptyset $ \ les arguments de [B.97]  s'appliquent et on a les isomorphismes suivants :
$$ \mathbb{H}_{c/f}^{n+1}(\mathcal{U}' -  \overline{\mathcal{U}} , \mathcal{E}^{\cdot}(1) , \delta_{0} )   \cong   \mathbb{H}_{c/f}^{n+1}(Y' ,\mathcal{E}^{\cdot}(1) , \delta_{0} )  \cong  H_{c}^{n+1}(Y',\mathbb{C}) \oplus  H_{c}^{n}(Y',\mathbb{C}) \frac{df}{f}   . $$
On est d'ailleurs ici dans une situation plus simple que celle de  [B.97] puisque  la fonction  $ f_{\mathbb{C}}$ \ est non singuli\`ere sur \  $\mathcal{U}' $ \ .

\noindent On peut donc trouver (voir (5) ci-dessus ) des formes  $\mathcal{C}^{\infty} , \eta$ \ , $ \omega $ \ sur  $\mathcal{U}' $ \ nulles au voisinage de  $ \overline{\mathcal{U}}$  et \`a supports  f-propre , qui sont d-ferm\'ees et de degr\'es respectifs  $n$  et  $n+1$ , ainsi qu'une forme semi-m\'eromorphe  \ $\alpha$ \ sur  $\mathcal{U}' $ \ , de degr\'e \ $ n $ \ , nulle au voisinage de  $ \overline{\mathcal{U}}$ \  et \`a supports  f-propre , v\'erifiant  sur  $\mathcal{U}' $   :
$$  W = \frac{df}{f}\wedge \eta  +  \omega  +  d\alpha  \  \  .$$ 
Alors  $ \widetilde{var}(e)$ \ est repr\'esent\'ee par
$$ \check{v} := (v_{1},\ldots ,v_{k-1}, v_{k} + \eta )  $$
o\`u \  $\eta$ \ est prolong\'ee par  $0$  \`a  $X_{\mathbb{C}}$ \  , c'est \`a dire que l'on a  $ \delta_{0}\check{v} = 0 $  \ ainsi que  $ v_{k} + \eta \vert_{F}  =  \widetilde{var}(e) $ \  et \  $ Supp(\check{v}) \cap \mathcal{U} = \emptyset  $  (ce qui montre que   $ v_{k} + \eta \vert_{F} $ \ donne bien une classe dans  $ H^{n}(F,F\cap \mathcal{U})_{=1} $ ) .

\noindent On a donc maintenant 
$$ \mathcal{I} ( \widetilde{var_{c}}(\Gamma(A)), \bar{e}) = \langle \Gamma(A),\overline{\widetilde{var}(e)} \rangle  = \frac{1}{(2i\pi )^{n}} \int_{\Gamma(A)} v_{k} + \eta  \  \   . $$
Maintenant le calcul suivant (analogue a celui de [B.97]  p.20-21 )  donne l'\'egalit\'e 
$$  Res ( \lambda = 0 ,\frac{1}{\Gamma(\lambda)}\int_{A}  f^{\lambda} \ \frac{df}{f} \wedge \sigma w_{k} )  =  -  Res ( \lambda = 0 ,\frac{1}{\Gamma(\lambda)}\int_{A}  f^{\lambda} \  \frac{df}{f} \wedge (v_{k} + \eta ) ) \  \ . $$
En effet, on a  sur  $ X_{\mathbb{R}} $ \  , pour $ \Re \lambda \gg 1  $
$$ d(f^{\lambda} \sigma (w_{k} + \frac{df}{f} \wedge (\rho \xi_{k}))) =    \lambda \frac{df}{f} f^{\lambda} \wedge \sigma w_{k} +  f^{\lambda}d \sigma \wedge (w_{k} + \frac{df}{f} \wedge (\rho \xi_{k}) )+ f^{\lambda} \frac{df}{f} \wedge v_{k} $$
car on a  $ \sigma v_{k} = v_{k}  $ \ .

\noindent En utilisant la formule  (5)  on obtient  :
$$  \begin{array}{l} Res ( \lambda = 0 ,\frac{1}{\Gamma(\lambda)}\int_{A}  f^{\lambda} \  \frac{df}{f} \wedge w_{k} )   \  + \\  Res ( \lambda = 0 ,\frac{1}{\Gamma(\lambda)}\int_{A}  f^{\lambda} \  \frac{df}{f} \wedge \eta )  \  + \\  Res ( \lambda = 0 ,\frac{1}{\Gamma(\lambda)}\int_{A}  f^{\lambda} \  \omega )  \  +   \\ Res ( \lambda = 0 ,\frac{1}{\Gamma(\lambda)}\int_{A}  f^{\lambda} \ d \alpha )  \  +  \\  Res ( \lambda = 0 ,\frac{1}{\Gamma(\lambda)}\int_{A}  f^{\lambda} \  \frac{df}{f} \wedge v_{k} )   =  0
\end{array} $$
Mais comme  $ \omega $ \ est  $ \mathcal{C}^{\infty} $ \ sur  $ X_{\mathbb{R}} $ \  et identiquement nulle pr\`es de l'origine $\frac{1}{\Gamma(\lambda)} \int_{A} f^{\lambda}  \omega  $ \  est une fonction enti\`ere . Comme  $ \alpha $  est semi-m\'eromorphe \`a p\^oles dans  $ f^{-1}(0) $ et est identiquement nulle pr\`es de l'origine 
 $$\frac{1}{\Gamma(\lambda)} \int_{A} f^{\lambda} d \alpha $$ 
 est \'egalement une fonction enti\`ere de  $ \lambda $ . Il nous reste alors l'\'egalit\'e d\'esir\'ee .

\noindent On ach\`eve la preuve en consid\'erant les fonctions  $ s \rightarrow \int_{f_{\mathbb{R}}=s \cap A } v_{j}  \  \  \forall j \in [1,k-1] $  \ et  \ $ s \rightarrow  \int_{f_{\mathbb{R}}=s \cap A } v_{k} + \eta  $ \ et en raisonnant comme dans la preuve du th\'eor\`eme 6.1 b)  de [B.M.02]  $   \qquad \blacksquare $

\newpage

\section {\bf Le cas $\partial A \subset \lbrace 0 \rbrace $.}

\bigskip

\noindent Nous allons maintenant consid\'erer des  \  $ A \in  H^{0}(X_{\mathbb{R}}- f_{\mathbb{R}}^{-1}(0) ,\mathbb{C} )  $ \  dont le bord  $ \partial A $ \ est concentr\'e \`a l'origine ; ceci revient \`a dire que  $ A $ \ est dans l'image de la restriction (injective) 
$$  i :  H^{0}(X_{\mathbb{R}}- \lbrace 0 \rbrace,\mathbb{C}) \longrightarrow  H^{0}(X_{\mathbb{R}}- f_{\mathbb{R}}^{-1}(0) ,\mathbb{C} ) \ \  . $$
Nous allons commencer par "\'etendre" l'application  $ \Gamma $ \ en une application  \ 
$$ \widehat{\Gamma} :   H^{0}(X_{\mathbb{R}}- \lbrace 0 \rbrace,\mathbb{C}) \longrightarrow  H_{c}^{n}(F)  $$
rendant commutatif le diagramme  suivant , not\'e \  $ \mathcal{D} $  \ dans la suite :


\bigskip

$$ \xymatrix{0 \ar[r] &  H^{0}(X_{\mathbb{R}}- \lbrace 0 \rbrace) \ar[d]^{\widehat{\Gamma}} \ar[r]^-{i} & H^{0}(X_{\mathbb{R}}- f_{\mathbb{R}}^{-1}(0)  ) \ar[d]^{\Gamma} \ar[r]  & H^{0}( f_{\mathbb{R}}^{-1}(0) -  \lbrace 0 \rbrace) \ar@{-->}[d]^{\cap \partial B(0,\varepsilon'')} \\
H^{n}(F\cap \mathcal{U}) \ar[r] &  H^{n}_{c}(F) \ar[r]^{can_{c}} & H^{n}_{c}(F,F\cap \mathcal{U}) \ar[r]^{\partial_{c}} &  H^{n+1}_{c}(F\cap \mathcal{U}) } $$

\bigskip

\noindent Nous montrerons ensuite que , quand \ $ \partial A \subset \lbrace 0 \rbrace $ \ , les p\^oles simples aux entiers n\'egatifs  du prolongement m\'eromorphe de la distribution  $ \int_{A} f^{\lambda} \Box $ \ sont contr\^ol\'es par la classe de cohomologie  $ \widehat{\Gamma}(A) $ \ (en fait par sa composante sur  \ $ H_{c}^{n}(F)_{=1} ) $ \ .

\bigskip

\bigskip

\noindent {\bf Remarques}

\begin{itemize}
\item {1)}   \  Comme  $ F\cap \mathcal{U} $ \ est homotopiquement \'equivalent \`a  $ \Lambda : = f_{\mathbb{R}}^{-1}(0) \cap \partial B(0,\varepsilon'') $ \ qui est lisse et compacte (orient\'ee) de dimension  $ n-1 $ \ , on a  $ H^{n}( F\cap \mathcal{U}) = 0  $ \ et donc aussi  \ $ H_{c}^{n}( F\cap \mathcal{U}) = 0  $ \ par dualit\'e de Poincar\'e . On a donc unicit\'e d'une telle application  \ $ \widehat{\Gamma} $ \  .

\item{2)}   \ Une fa{\c c}on "peu explicite" (c'est \`a dire en restant au niveau cohomologique sans exhiber un  n-cycle compact \ $ \widehat{\Gamma}(A) $ \  pour chaque \  $A$ \ donn\'e ) de construire  \ $ \widehat{\Gamma}$ \  est de voir que  $ \Gamma $ \  commute \`a l'intersection avec  $  \partial B(0,\varepsilon'') $ \ qui d\'efinit une application lin\'eaire  $ H^{0}(f_{\mathbb{R}}^{-1}(0) - \lbrace 0 \rbrace , \mathbb{C}) \rightarrow H_{c}^{n+1}(F\cap \mathcal{U} ) $ \ . On utilise ici l'isomorphisme
$$ H^{1}_{f_{\mathbb{R}}^{-1}(0)}(X_{\mathbb{R}}- \lbrace 0 \rbrace,\mathbb{C}) \cong  H^{0}( f_{\mathbb{R}}^{-1}(0) -  \lbrace 0 \rbrace ,\mathbb{C}) $$
qui est cons\'equence de la lissit\'e de \ $  f_{\mathbb{R}}^{-1}(0) -  \lbrace 0 \rbrace $ \ dans \ $ X_{\mathbb{R}}- \lbrace 0 \rbrace $ \ via la d\'eg\'en\'erescence de la suite spectrale de cohomologie \`a support  . L\`a encore,il s'agit d'un s\'erieux coup de chapeau \`a Jean Leray !
Cette compatibilit\'e est facile \`a voir sur notre construction de  $ \Gamma $ \ . En effet,pour  $ A_{\alpha} \subset \lbrace f_{\mathbb{R}} > 0 \rbrace $ c'est \'evident ; pour  \  $ A_{\alpha} \subset \lbrace f_{\mathbb{R}} < 0 \rbrace $ \ ceci utilise la compatibilit\'e entre la  trivialisation  $\Phi $ \ de  $ f \vert_{\mathcal{U}'} $ \ et la trivialisation \ $ \Psi $ que l'on suit le long du demi-cercle  \ $ -s_{0}.e^{i\theta} ,$ \ pour \ $\theta \in [-\pi,0] $ \ pour amener \ $ \Gamma(A)^{-} $ dans la fibre de Milnor  \ . Alors le diagramme \  $ \mathcal{D} $ \ permet facilement de construire  l'application \ $ \widehat{\Gamma} $ \ .

\end{itemize}

\noindent Nous allons expliciter compl\`etement  le cycle  \ $ \widehat{\Gamma}(A) $ \  (c'est \`a dire construire un n-cycle compact de  $F$ )  quand  $ \partial A \subset \lbrace 0 \rbrace $ \ et prouver le 

\bigskip

\noindent{\bf Theoreme 2}

\noindent \textit{Il existe une unique application lin\'eaire  $ \widehat{\Gamma} $ \  rendant commutatif le diagramme \  $ \mathcal{D} $  \  ci-dessus . }

\noindent \textit{Pour  $ A \in H^{0}(X_{\mathbb{R}}- \lbrace 0 \rbrace, \mathbb{C}) $ \ et \ $ e \in H^{n}(F)_{=1} $ \ on a  :
\begin{equation}
  (-2i\pi)^{n}\  \mathcal{I} ( \widehat{\Gamma}(A) , \bar{e} )  \ =  \  Res (\lambda = 0 , \int_{A} (f/s_{0})^{\lambda}  \frac{df}{f} \wedge w_{k} ) 
\end{equation}
o\`u \  $ w_{1},\ldots,w_{k} $ \ repr\'esentent   \ $ e $ \ dans  $ H^{n}(F)_{=1} $  , c'est \`a dire v\'erifient la condition  (2)   du th\'eor\`eme 1  avec  $ u = 0    \qquad  \  \qquad  \  \Box $ }

\bigskip

\noindent {\bf Remarques}

\begin{itemize}
\item {1)}   \  Comme on a  \ $\partial A \subset \lbrace 0 \rbrace $ \ les parties polaires aux entiers n\'egatifs du prolongement m\'eromorphe de \  $ \int_{A} f^{\lambda} \Box $ \ sont toutes des distributions \`a support l'origine . C'est la raison pour laquelle on peut omettre dans la formule du th\'eor\`eme 2  de faire appara\^itre une fonction   \ $\rho \in \mathcal{C}^{\infty}(X_{\mathbb{R}}) $ \ valant identiquement   1  au voisinage de l'origine .

\item{2)}  \ Dans le cas o\`u   \ $\partial A \not\subset \lbrace 0 \rbrace $ \  il est facile de voir que les r\'esidus aux entiers n\'egatifs du prolongement m\'eromorphe de \  $ \int_{A} f^{\lambda}\  \Box $ \ sont des distributions non nulles le long du bord de \ $ A $ \ . Ces p\^oles (simples en dehors de $0$ )  ne sont pas li\'es \`a la singularit\'e de  $ f_{\mathbb{R}} $ \ et existent toujours  (et sont toujours non nuls )
le long de  \ $ \partial A - \lbrace 0 \rbrace $ \  . Dans cette situation les th\'eor\`emes 1  et  1 bis  d\'ecrivent compl\`etement les parties polaires "int\'eressantes" du prolongement  m\'eromorphe de \  $ \int_{A} f^{\lambda} \  \Box $ \ . Par contre, si  \ $\partial A \subset \lbrace 0 \rbrace $ \ , le th\'eor\`eme 2  est n\'ecessaire pour compl\'eter la description des p\^oles "int\'eressants".

\item{3)} La commutativit\'e du diagramme ci-dessus implique l'\'egalit\'e \ $$ can_{c}(\widehat{\Gamma}(A) = \Gamma(A) $$ \ et donc aussi \ $$ var_{c}(\Gamma(A)) = (T - 1)(\widehat{\Gamma}(A)) . $$
\end{itemize}

\bigskip

\noindent{\bf Corollaire }

\noindent\textit{ Dans la situation du Th\'eor\`eme 2 , le prolongement m\'eromorphe de \ $ \int_{A} f^{\lambda} \ \Box $ \ n'a pas de p\^oles aux entiers n\'egatifs si et seulement si   \ $ \widehat{\Gamma}(A)_{=1} =  0  $ \ dans \ $ H_{c}^{n}(F)_{=1} .$ }
\noindent \textit{De plus, le prolongement  m\'eromorphe de \ $ \int_{A} f^{\lambda} \ \Box $ \ n'a pas de p\^oles du tout si et seulement si \ $ \widehat{\Gamma}(A)  =  0  $ \ dans \ $ H_{c}^{n}(F)  \  \  \  \  \  \Box  $ }

\bigskip

\bigskip

\noindent Preuve du th\'eor\`eme 2 : 

\noindent Compte tenu des outils d\'ej\`a mis en place (y compris dans l'Annexe), cette preuve suit pas \`a pas celle de [B.02] . Aussi nous contenterons-nous de l'esquisser pour la commodit\'e du lecteur , renvoyant \`a  [B.02]  pour plus de d\'etails .

\begin{itemize}
\item{a)} Construction de $ \widehat{\Gamma}(A) $ :

\noindent Pour \ $ A \in H^{0}(X_{\mathbb{R}}- \lbrace 0 \rbrace,\mathbb{C}) $ \ , \ $ \Gamma(A)^{+} $ \ et \  $ \Gamma(A)^{-} $ \ sont deux $n-$cha\^ines compactes orient\'ees de  \ $X_{\mathbb{R}} $ \ \`a bords dans  \ $ \mathcal{U} $ \ . Leurs bords \ $ \partial  \Gamma(A)^{+} $ \ et \ $ \partial  \Gamma(A)^{-} $ \
 sont deux \ $(n-1)$-cycles compacts (orient\'es) de \  $\mathcal{U} $ \ qui sont homologues dans \ $\mathcal{U} $ \ . L'homologie est donn\'ee par la \ $n-$cha\^ine compacte de \ $\mathcal{U} $ \ :
$$  \Delta =  A \cap f_{\mathbb{R}}^{-1}([-s_{0},s_{0}]) \cap \partial B(0,\varepsilon'')  $$
dont le bord est 
$$   \partial  \Gamma(A)^{+} -  \partial  \Gamma(A)^{-} + (\sum_{\alpha}a_{\alpha}.\bar{A_{\alpha}}) \cap  f_{\mathbb{R}}^{-1}(0) \cap \partial B(0,\varepsilon'') $$
mais l'hypoth\`ese \ $\partial A \subset \lbrace 0 \rbrace $ \  donne 
$$( \sum_{\alpha}a_{\alpha}. \bar{A_{\alpha}}) \cap  f_{\mathbb{R}}^{-1}(0) \cap \partial B(0,\varepsilon'')  = \emptyset  \ \ .$$
Soit \ $ \Delta_{0} $ \ la  $n-$cha\^ine compacte (orient\'ee) de \ $ \mathcal{U} $ \ qui est l'image r\'eciproque par la restriction de  \ $ f_{\mathbb{R}} $ \ \`a \ $ \partial B(0,\varepsilon'') $ \ 
 du chemin \ $\gamma_{0}$ \ obtenu en d\'eformant le segment  \ $ [-s_{0},s_{0}] $ \  par un  demi-cercle contournant l'origine par en-dessous ( on utilise ici la trivialisation \ $\Phi$ \ ) :


\noindent Alors le  $n-$cycle compact  \ $ \Gamma(A)_{0} : = \Gamma(A)  -  \Delta + \Delta_{0} $ \ de \
$$ f_{\mathbb{C}}^{-1}( (\bar{D} -  0 ) \cap \lbrace \Im \leq 0 \rbrace)  $$
 se d\'eforme, en utilisant la trivialisation \ $ \Psi $ cette fois , en suivant la d\'eformation suivante du chemin \ $\gamma_{0}$ \ :

\begin{center}
\includegraphics[width=100mm]{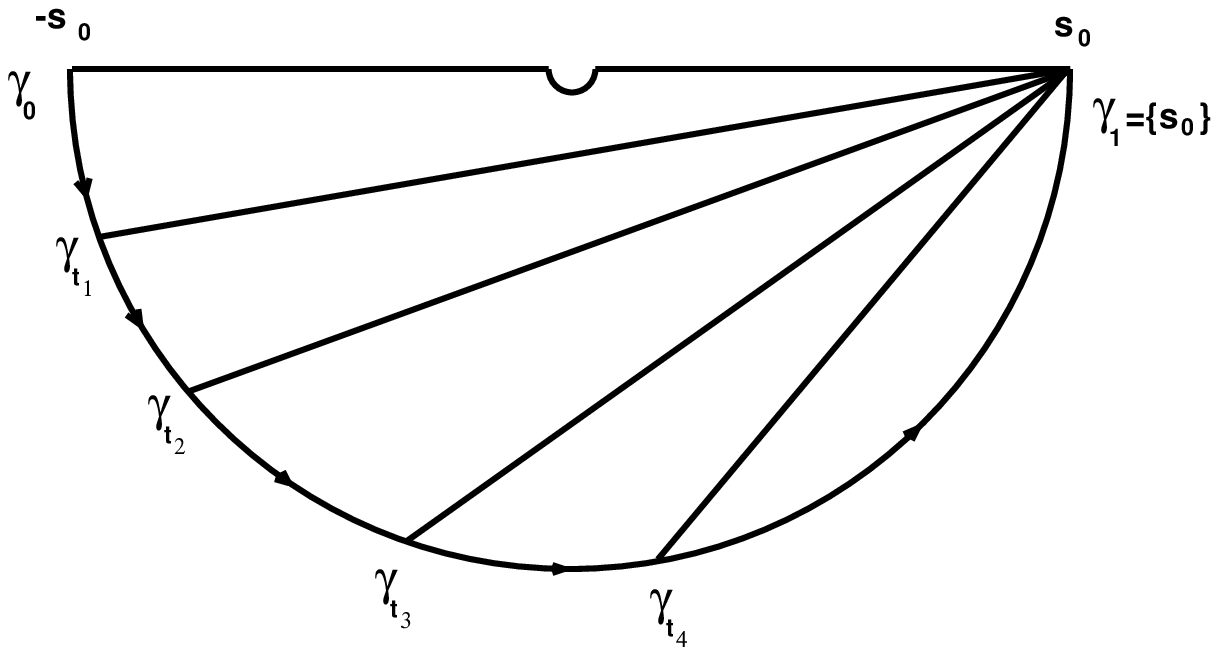}
\end{center}

\noindent o\`u \ $\gamma_{1}$ \ est le chemin constant \'egal \`a \ $ s_{0} $ \ et \ $ \gamma_{t}(0) = - s_{0} e^{i\pi t} $ \ . Alors on a 
$$ \Gamma(A)_{t} =  \Gamma(A)^{+}  - T^{\frac{t}{2}}\Gamma(A)^{-}  + \Delta_{t}  $$
La notation \ $T^{\frac{t}{2}} $ \ signifie que l'on a suivi le $\frac{1}{2}-$cercle \ $ -s_{0} e^{i\pi \theta}$ \ pour \ $ \theta \in [0,t] $ \ dans la trivialisation \ $\Psi $ \  et \ $ \Delta_{t} $ \ d\'esigne la d\'eform\'ee de \ $\Delta_{0} $ \ en suivant la d\'eformation de  \ $\gamma_{0} $ \ \`a \ $ \gamma_{t} $ \ .

\noindent On pose alors \ $ \widehat{\Gamma}(A) = \Gamma(A)_{1}  \subset  F = f_{\mathbb{C}}^{-1}(s_{0}) $ \ .

\noindent Comme la trivialisation \ $\Psi $ \ est compatible avec la trivialisation \  $\Phi $ \ de \ $ f\vert_{\mathcal{U}} $ \ , on a 
$$  can_{c}(\widehat{\Gamma}(A)) = \Gamma(A) $$
dans \ $ H_{c}^{n}(F,F\cap \mathcal{U} ) $ \ et pour achever la d\'emonstration du th\'eor\`eme 2 il nous reste seulement \`a prouver la formule (6) .

\bigskip

\item{b)} Preuve de (6) :

\noindent Consid\'erons donc \ $ e \in H^{n}(F)_{=1} $ \ et \ $ w_{1},\ldots,w_{k} $ \ des  $n-$formes semi-m\'eromorphes v\'erifiant la condition (2) du th\'eor\`eme 1  avec  $u=0 $ \ ( donc  $ w_{k}$ \  induit  $ e $ \ dans \ $  H^{n}(F) $ \ ) . Fixons la d\'etermination du logarithme sur  $\mathbb{C} - i\mathbb{R}^{+} $  de facon que  $ arg(z) \in ]-\frac{3\pi}{2},\frac{\pi}{2}[ $ \ et posons  :
$$\Omega : = \sum_{j=0}^{k-1} \frac{(-1)^{j}}{j!} \lbrack Log(f/s_{0}) \rbrack^{j} \  w_{k-j} $$
La n-forme  $\Omega $ \ est semi-m\'eromorphe sur  \ $ X_{\mathbb{C}} - f^{-1}(D\cap i\mathbb{R}^{+}) $ \ et on a \ $ d\Omega = 0 $ \ et \ $ \Omega \vert_{F} = w_{k} \vert_{F} = e $ \ .

\noindent Comme les $n-$cycles compacts \ $ \Gamma(A)_{0} $ \ et \ $ \Gamma(A)_{1} = \widehat{\Gamma}(A) $ \ sont homologues dans \ $ X_{\mathbb{C}} - f^{-1}(D\cap i\mathbb{R}^{+}) $ \ on a  :
$$(2i\pi)^{n} \  \mathcal{I} ( \widehat{\Gamma}(A) , \bar{e} ) =  \int_{\widehat{\Gamma}(A)} e  = \int_{\gamma(A)_{0}} \Omega  \  . $$
Maintenant les m\^emes arguments que [B.02] p.8  permettent de voir que l'on a  :
$$ \int_{\gamma(A)_{0}} \Omega = \int_{\partial [A\cap( f_{\mathbb{R}})^{-1}[-s_{0},s_{0}]\cap \overline{B(0,\varepsilon'')}]}  \Omega  \  \  .$$
Il reste alors \`a montrer que cette derni\`ere int\'egrale coincide bien avec 
$$ Res(\lambda = 0 , \int_{A\cap \overline{ B(0,\varepsilon'')}}  f^{\lambda} \ \frac{df}{f} \wedge w_{k} \ ) $$
ce qui donnera le r\'esultat gr\^ace au lemme 2  de [B.02] .

\noindent Comme on a  \ $ d( (f/s_{0})^{\lambda} \  \Omega ) = \lambda  \frac{df}{f} \wedge  (f/s_{0})^{\lambda} \  \Omega $ \ , la formule de Stokes donne  :
$$ \int_{A\cap( f_{\mathbb{R}})^{-1}[-s_{0},s_{0}]\cap \overline{B(0,\varepsilon'')}]} \lambda  \ \frac{df}{f} \wedge  (f/s_{0})^{\lambda} \  \Omega  =  \int_{\partial [A\cap( f_{\mathbb{R}})^{-1}[-s_{0},s_{0}]\cap \overline{B(0,\varepsilon'')}]}  (f/s_{0})^{\lambda}  \ \Omega \ $$
ce qui donne , en \ $ \lambda = 0 $ \ :
$$  Res(\lambda = 0 , \int_{A\cap \overline{ B(0,\varepsilon'')}}  f^{\lambda} \ \frac{df}{f} \wedge \Omega \ )  =  \int_{\partial [A\cap( f_{\mathbb{R}})^{-1}[-s_{0},s_{0}]\cap \overline{B(0,\varepsilon'')}]}  \Omega  \  \  .$$
Revenons \`a la d\'efinition de \ $ \Omega $ \ : le membre de gauche de l'\'egalit\'e ci-dessus vaut donc 
$$Res(\lambda = 0,  \sum_{j=0}^{k-1} \frac{(-1)^{j}}{j!} \int_{A\cap \overline{B(0,\varepsilon'')}} (f/s_{0})^{\lambda}\lbrack Log(f/s_{0}) \rbrack^{j}  \ \frac{df}{f} \wedge w_{k-j}  )     \ .     $$
Mais pour \ $ j\geq 1$ \ on a 
$$ \frac{d}{d\lambda} \lbrack \int_{A\cap \overline{ B(0,\varepsilon'')}} (f/s_{0})^{\lambda}\lbrack Log(f/s_{0}) \rbrack^{j-1} \rho  \ \frac{df}{f} \wedge w_{k-j} \rbrack =  $$
 $$ \int_{A\cap \overline{ B(0,\varepsilon'')} }(f/s_{0})^{\lambda}\lbrack Log(f/s_{0}) \rbrack^{j} \rho  \ \frac{df}{f} \wedge w_{k-j}  $$
et la d\'eriv\'ee d'une fonction m\'eromorphe n'a jamais de r\'esidu . Il reste donc seulement le terme en \ $ j = 0 $ \ et ceci ach\`eve la d\'emonstration  \  \  $\blacksquare $

\end{itemize}


\newpage

\section{\bf Annexe}

\noindent Soit  \ $ X $ \ une vari\'et\'e \ $\mathcal{C}^{\infty} $ \ paracompacte et soit \ $ \mathcal{U} $ \ un ouvert de   \ $ X $ \ . On suppose qu'il existe deux voisinages ouverts  de \ $ \overline{ \mathcal{U}} \ ,  \ \mathcal{U}''\ $ \ et \ $ \mathcal{U}' $ \ v\'erifiant   $\overline{ \mathcal{U}'' } \subset  \mathcal{U}' $ \ et tels que les inclusions  \ $   \mathcal{U}\subset  \mathcal{U}'' \subset  \mathcal{U}' $ \ soient des \'equivalences d'homotopie .
\noindent Dans cette situation nous d\'efinirons , pour \ $ p \in \mathbb{N} $ 
$$  H_{c}^{p}(X,\mathcal{U}) : = \genfrac{}{}{}{}{\lbrace \varphi \in \mathcal{C}_{c}^{\infty}(X)^{p} / Supp(d\varphi )\subset \mathcal{U} \rbrace}{ \mathcal{C}_{c}^{\infty}(\mathcal{U})^{p} + d\mathcal{C}_{c}^{\infty}(X)^{p-1}}  $$

\noindent{\bf Lemme 1}

\noindent \textit{ On a une suite exacte longue de cohomologie  :
$$ \cdots  \xrightarrow{\gamma} H_{c}^{n}(X) \xrightarrow{\alpha}H_{c}^{n}(X,\mathcal{U})\xrightarrow{\beta}H_{c}^{n}(\mathcal{U})\xrightarrow{\gamma} H_{c}^{n+1}(X)\xrightarrow{\alpha} \cdots $$
o\`u \ $ \alpha $ \ est l'application canonique ( $ d\varphi = 0 $ \  implique \ $ Supp (d\varphi) \subset \mathcal{U} $ \ ) , \ $ \beta[\varphi] = [d\varphi] $ \ et \ $ \gamma $ \ est le prolongement par \ $0$ .}

\bigskip

\noindent La preuve est \'el\'ementaire .

\noindent D\'efinissons \'egalement pour \ $ p \in \mathbb{N} $ 
$$  H^{p}(X,\mathcal{U}) : = \genfrac{}{}{}{}{\lbrace \varphi \in \mathcal{C}^{\infty}(X)^{p} / d\varphi = 0  \ et \ Supp(\varphi )\cap \mathcal{U} = \emptyset \rbrace}{ d(\mathcal{C}^{\infty}(X)^{p-1}\cap \lbrace Supp\cap \mathcal{U} = \emptyset \rbrace)}  $$

\bigskip

\noindent{\bf Lemme 2}

\noindent\textit{ On a une suite exacte longue de cohomologie  :
$$ \cdots \xrightarrow{\hat{\alpha}} H^{m-1}(X) \xrightarrow{\hat{\gamma}}H^{m-1}(\mathcal{U})\xrightarrow{\hat{\beta}}H^{m}(X,\mathcal{U}) \xrightarrow{\hat{\alpha}} H^{m}(X)\xrightarrow{\hat{\gamma}}\cdots $$
o\`u \ $\hat{ \alpha} $ \ est l'application canonique ( oublie de la condition de support) , o\`u \ $ \hat{\gamma}$ \ est la restriction et o\`u \ $\hat{\beta} $ \ est d\'efinie  ci-dessous  .} 
\bigskip

\noindent Fixons une fonction \ $\rho \in \mathcal{C}^{\infty}(X) $ \ identiquement \'egale a  1  sur \ $ X - \overline{\mathcal{U}''} $ \ et identiquement nulle sur  \ $\mathcal{U} $ . Si \ $ [\varphi] \in H^{m-1}(\mathcal{U})$ \ , on utilise l'isomorphisme induit par la restriction \  $H^{m-1}(\mathcal{U}')\rightarrow  H^{m-1}(\mathcal{U}) $ \ (gr\^ace \`a notre hypoth\`ese d'\'equivalence d'homotopie) pour trouver une \ $(m-1)-$forme \ $ \mathcal{C}^{\infty} \ \varphi '$ \ sur \ $\mathcal{U}' $ \ d-ferm\'ee et induisant sur  \ $\mathcal{U} $ \ la classe \ $[\varphi]$ \ . Posons alors  :
$$  \hat{\beta}[\varphi] = - [ d\rho \wedge \varphi' ]  \  \ .$$
La  $m-$forme \  \ $ \mathcal{C}^{\infty} \  d\rho \wedge \varphi'  $ \ \ est bien d-ferm\'ee et son support ne rencontre pas \  $\mathcal{U} .$ \ Comme la classe de cohomologie de  \  $\mathcal{U}' $ \ d\'efinie par \ $\varphi' $ \ est bien d\'etermin\'ee, on voit facilement qu'ajouter  \ $d\psi' $ \ \`a \ $ \varphi' $ \ ne fait qu'ajouter le terme \ $ d(d\rho \wedge \psi'))$ \  \`a \ $  -  d\rho \wedge \varphi' $ \ ce qui ne change pas la classe dans \ $H^{m}(X,\mathcal{U}) $ \ puisque le support de \ $d\rho $ \ ne rencontre pas \ $ \mathcal{U}  . $

\bigskip

\noindent Preuve  :

\noindent Nous allons nous contenter de v\'erifier l'inclusion \ $ Ker( \hat{\alpha}) \subset Im(\hat{\beta}) $  \ qui est le seul point non trivial .

\noindent Soit \ $ \varphi \in \mathcal{C}^{\infty}(X)^{m} $ \ v\'erifiant \ $Supp(\varphi)\cap \mathcal{U} = \emptyset $ \ et d-exacte sur  X  ( donc \ $ \varphi \in  Ker( \hat{\alpha}) $ \ ) . Soit  alors \ $\psi \in \mathcal{C}^{\infty}(X)^{m-1} $ \ telle que \ $ d\psi = \varphi $ . Nous voulons v\'erifier que \ $ [\varphi] = \hat{\beta}(\psi \vert_{\mathcal{U}}) $ ( on remarquera que \ $(\psi \vert_{\mathcal{U}} $ \ est d-ferm\'ee ) . 
Par d\'efinition , on a \ $  \hat{\beta}(\psi \vert_{\mathcal{U}}) = -d\rho \wedge \psi' $ \ ou \ $ \psi' \in \mathcal{C}^{\infty}(\mathcal{U}')^{m-1} $ \ est d-ferm\'ee  et v\'erifie , quitte \`a changer     le repr\'esentant de la classe \ $ [\psi] $ , \ $ \psi' \vert_\mathcal{U} =  \psi \vert_\mathcal{U} $ \ . La forme \ $ \chi : = \psi - (1- \rho) \psi' $ \ est \ $  \mathcal{C}^{\infty} $ \ sur \ $X$ \ de degre \ $m-1 $ \ , v\'erifie \ $ Supp(d\chi) \cap \mathcal{U} = \emptyset $ \ et on aura dans \ $ H^{m}(X,\mathcal{U}) $ \ :
$$ [0] = [d\chi] = [d\psi] - [-d\rho \wedge \psi'] = [\varphi] - [\hat{\beta}(\psi \vert_{\mathcal{U}})]  \ \  \qquad  \  \  \qquad  \  \  \blacksquare$$

\bigskip

\noindent Supposons maintenant que \ $X$  soit orient\'ee et de dimension \ $n$ \ . On a pour chaque \ $ p \in \mathbb{N} $ \ un accouplement  sesquilin\'eaire  
$$\langle , \rangle : \  H^{p}(X,\mathcal{U}) \times H_{c}^{n-p}(X,\mathcal{U}) \longrightarrow \mathbb{C}  $$
d\'efini par
$$  \langle  [\psi] , [\varphi] \rangle  = \frac{1}{(2i\pi)^{n}} \int_{X}  \psi \wedge \overline{\varphi}  $$
o\`u  \ $ \psi  \in \mathcal{C}^{\infty}(X)^{p} $ \ v\'erifie \ $ Supp(\psi) \cap \mathcal{U} = \emptyset  $ \ et \ $ d\psi = 0 $ \ et o\`u \ $  \varphi \in \mathcal{C}^{\infty}(X)^{n-p} $ \ v\'erifie \ $ Supp(d\varphi) \subset \mathcal{U} $ \ . 

\bigskip


\noindent Preuve :

\noindent Soient \ $ \eta \in \mathcal{C}^{\infty}(X)^{p-1} $ \ telle que \ $ Supp(\eta) \cap \mathcal{U} = \emptyset , \  \xi \in  \mathcal{C}_{c}^{\infty}(\mathcal{U})^{n-p} \ $ \ et \ $ \zeta \in  \mathcal{C}_{c}^{\infty}(X)^{n-p-1} \ $ \ . Il s'agit de v\'erifier l'\'egalit\'e  :
$$ \int_{X}  \psi \wedge \overline{\varphi} = \int_{X} (\psi + d\eta) \wedge \overline{(\phi + \xi  + \ d\zeta )}  $$
ce qui revient \`a montrer que 
$$ \int_{X} \psi \wedge \overline{(\xi + d\zeta)}  =  0  \ \ \  \  et \  \   \  \  \int_{X} d\eta \wedge \overline{(\varphi + \xi + d\zeta )}  = 0   \  .$$
Comme \ $ Supp(\psi) \cap \mathcal{U} = \emptyset  , \  Supp(\eta) \cap \mathcal{U} = \emptyset  $ \ et \ $ \xi \in \mathcal{C}_{c}^{\infty}(\mathcal{U})^{n-p} $ \ on a 
$$ \int_{X} \psi \wedge \overline{\xi}  = 0 = \int_{X} d\eta \wedge \overline{\xi}  \  \ . $$
De plus la formule de Stokes donne imm\'ediatement la nullit\'e de \ $ \int_{X} d\eta \wedge \overline{d\zeta}  .$ \  Elle donne aussi la nullit\'e de \ $ \int_{X} \psi \wedge \overline{d\zeta} $ \ puisque \ $\psi $ \ est d-ferm\'ee .

\noindent Enfin \ $ \int_{X} d\eta \wedge \overline{\varphi} = \pm  \int_{X} \eta \wedge \overline{d\varphi} = 0 $ \ car \ $ Supp(d\varphi) \subset \mathcal{U} \  \   \  \  \  \  \  \  \  \  \blacksquare $

\bigskip 

\noindent{\bf Proposition}

\noindent \textit{Sous nos hypoth\`eses , l'accouplement \ $ \langle , \rangle$ \ est compatible aux "dualit\'es"de Poincar\'e  (hermitiennes) sur \ $ X $ \ et \ $\mathcal{U} $ \ . Si l'on suppose de plus les cohomologies de  \ $ X $ \ et \ $\mathcal{U} $ \ sont de dimensions finies , alors l'accouplement \ $\langle , \rangle $ \ est une dualit\'e hermitienne $ \  \  \  \  \  \  \  \  \  \    \  \  \  \  \  \  \  \  \  \Box $ }

\bigskip

\noindent Preuve :

\noindent Notons par \ $  \mathcal{I}  $  \ les "dualit\'es" de Poincar\'e ( sur \  $ X $ \ et \ $\mathcal{U} $ \ ) . La premi\`ere compatibilit\'e de l'\'enonce signifie que l'on a l'\'egalit\'e 
$$\langle \psi, \alpha(\varphi) \rangle =  \mathcal{I} ( \hat{\alpha}(\psi) , \varphi ) $$
qui est imm\'ediate ; la seconde compatibilit\'e est donn\'ee par  l'\'egalit\'e
$$  (-1)^{deg(\psi)}  \langle \hat{\beta}(\psi) , \varphi \rangle = \mathcal{I}( \psi , \beta(\varphi))  $$
qui s'obtient de la facon suivante :  
$$ (2i\pi)^{n} \ \mathcal{I}( \psi , \beta(\varphi) ) = \int_{\mathcal{U}} \psi \wedge \overline{d\varphi}  =  \int_{\mathcal{U'}} \psi' \wedge \overline{d\varphi} = \int_{\mathcal{U'}} \psi' \wedge (1- \rho ) \overline{d\varphi} $$
car \ $\rho \equiv 0 $ \ sur \ $ Supp(d\varphi) \subset \mathcal{U} $ \ . On en d\'eduit :
$$  (2i\pi)^{n}\ \mathcal{I}( \psi , \beta(\varphi) ) =  (-1)^{deg(\psi)} \int_{\mathcal{U}'} - d\rho \wedge \psi' \wedge \overline{\varphi} =  (2i\pi)^{n}  (-1)^{deg(\psi)} \langle \hat{\beta}(\psi) , \varphi \rangle $$
ce qui prouve les compatibilit\'es d\'esir\'ees .

\bigskip

\noindent Prouvons maintenant la non-d\'eg\'en\'erescence de  \ $  \langle , \rangle$ . Il suffit de v\'erifier que si \ $ e \in H^{p}(X,\mathcal{U}) $ \ v\'erifie \ $ \forall \varepsilon \in H_{c}^{n-p}(X,\mathcal{U}) , \ \mathcal{I}(e,\varepsilon) = 0 $ \ alors on a \ $ e = 0 $ \ et que si \ $ \varepsilon \in H_{c}^{n-p}(X,\mathcal{U}) $ \ v\'erifie \ $\forall e \in H^{p}(X,\mathcal{U}) ,\  \mathcal{I}(e,\varepsilon) = 0 $ \ alors \ $ \varepsilon = 0 $ \ puisque les suites exactes longues de cohomologie donnent la finitude des espaces vectoriels (complexes) \ $ H^{*}(X,\mathcal{U}) \ $ \ et \ $ H_{c}^{*}(X,\mathcal{U}) \ $ \ sous nos hypoth\`eses  .

\noindent Comme il s'agit d'un simple exercice, traitons le premier point pour la commodit\'e du lecteur :  On a , par hypoth\`ese \ $ \forall \eta \in  H_{c}^{n-p}(X) , \langle e, \alpha(\eta) \rangle  = 0 $ \ et donc \ $ \hat{\alpha}(e) = 0 $ \ d'apr\`es la dualit\'e de Poincar\'e sur \ $ X $ \ et la premi\`ere compatibilit\'e ci-dessus . Il existe  donc  \ $ \zeta \in H^{p-1}(\mathcal{U}) $ \ tel que \ $ e = \hat{\beta}(\zeta) $ \ . Mais alors
$$ \langle \hat{\beta}(\zeta), \varepsilon \rangle  = (-1)^{p-1}\mathcal{I} (\zeta , \beta(\varepsilon))  = 0  \ . $$ Donc on a 
$$ \zeta \in (Im \beta )^{\bot} = (Ker \gamma)^{\bot} = Im \hat{\gamma} $$
car \ $\gamma$ \  et \ $ \hat{\gamma} $ \ sont adjoints via les dualit\'es de Poincar\'e sur \ $\mathcal{U}$  \ et \ $ X $ . On a donc \ $ \zeta = \hat{\gamma}(\xi) $ \ o\`u \ $\xi \in H^{p-1}(X) $ \ et donc \ $ e = \hat{\beta}(\hat{\gamma}(\xi)) = 0  \   \   \   \    \   \   \blacksquare$

\bigskip

\noindent{\bf Remarque}

\noindent Il est facile de voir que , comme dans la cohomologie de de Rham "standart" , on peut remplacer les formes \ $ \mathcal{C}^{\infty} $ \ par des courants (disons pour \ $ X $ \ orient\'ee, seul cas qui nous int\'er\`esse  ici) . Ceci permet alors de faire le lien avec un calcul via des cha\^ines orient\'ees  \`a bords dans \ $ \mathcal{U} $ \ . Ce point de vue qui est utilis\'e dans le texte pour d\'ecrire (rapidement) notre construction de l'application de variation , est une simple commodit\'e qui pourrait \^etre \'evit\'ee en utilisant directement le point de vue des formes \ $ \mathcal{C}^{\infty} $ \ . Bien sur dans ce cas l'op\'eration de troncature doit \^etre remplac\'ee par la multiplication par une fonction plateau qui est moins "parlante" (mais qui r\'eapparait finalement via le lemme 2 de [B.02] !) . Ceci conduirait \`a une d\'emonstration plus p\'enibles \`a suivre de nos th\'eor\`emes . C'est pourquoi nous avons donn\'e cette pr\'esentation de la variation , en nous attachant cependant \`a v\'erifier que son "adjoint "est bien donn\'e par le residu de J. Leray  de fa{\c c}on \`a pouvoir utiliser le calcul de la variation en cohomologie de de Rham donn\'e dans  [B.97] .

\newpage

\section {\bf R\'ef\'erences}

\bigskip

\begin{itemize}

\item{[A.70] Atiyah, M.F.  \textit{Resolution of singularities and division of distributions .} Comm.in pure and appl. Math. 23 (1970) p.145-150  .}

\item{[B.85]  Barlet, D.   \textit{La forme hermitienne canonique sur la cohomologie de la fibre de Milnor d'une hypersurface \`a singularit\'e isol\'ee .} Invent. Math. 81 (1985) p.115-153 .}

\item{[B.90]  Barlet, D.   \textit{La forme hermitienne canonique pour une singularit\'e presque isol\'ee .} Complex Analysis (K. Diederich eds) Vieweg,Wuppertal (1990) p.20-28  .}

\item{[B.91]  Barlet, D.   \textit{Emm\^elements de strates cons\'ecutives pour les cycles \'evanescents .} Ann. Scient. ENS 4-i\`eme s\'erie  24 (1991) p.401-506  .}

\item{[B.97]  Barlet, D.   \textit{La variation pour une hypersurface \`a singularit\'e isol\'ee relativement \`a la valeur propre  1 .} Revue de l'Inst. E. Cartan (Nancy) 15 (1997) p.1-29  .}

\item{[B.99]  Barlet, D.  \textit{Multiple poles at negative integers for \ $ \int_{A} f^{\lambda} \ \Box $ \ in the case of an almost isolated singularity .} Publ. RIMS,Kyoto Univ. 35 (1999)  p.571-584  .}

\item{[B.02]  Barlet, D.  \textit{ Real canonical cycle and asymptotics of oscillating integrals .} Pr\'epublication de l'Institut E. Cartan (Nancy) 2001 /35  \`a para\^itre au Nagoya Math.Journ.}


\item{[B.M.93]   Barlet, D.  et  Maire, H.-M.  \textit{Asymptotique des int\'egrales-fibres.} Ann. Inst.Fourier (Grenoble) 43 (1993) p.1267-1299  .}

\item{[B.M.00]  Barlet, D.  et  Maire, H.-M.   \textit{Poles of the current \ $ |f|^{2\lambda} $ \ over an isolated singularity .}  Intern. Journ. Math. 11 (2000) p.609-635  .}

\item{[B.M.02]  Barlet, D.  et  Maire, H.-M.   \textit{Non trivial simple poles at negative integers and mass concentration at singularity .} Math.Ann. 323 (2002) p.547-587  .}

\item{[H.64]  Hironaka, H.  \textit{Resolution of singularities of an algebraic variety over a field of characteristic zero .} Ann. of Math. 79 (1964) p.109-326 }

\item{[J.91]  Jeddi, A.   \textit{ Singularit\'e r\'eelle isol\'ee .} Ann.Inst.Fourier (Grenoble) 41 (1991) p.87-116 .}

\item{[J.02]  Jeddi, A.  \textit{ Preuve d'une conjecture de Palamodov .} Topology  41 (2002)  p.271-306 .}

\item{[JQ.70] Jeanquartier, P.  \textit{D\'eveloppements asymptotiques de la distribution de Dirac attach\'ee \`a une fonction analytique .} C.R.A.S. Paris 271 (1970) p.1159-1161 .}

\item{[L.59]  Leray, J.   \textit{Le probl\`eme de Cauchy III . } Bull.Soc. Math. France 87 (1959)  p.81-180  .}

\item{[Lj.65]  Lojasiewicz, S.   \textit{Ensembles semi-analytiques .} IHES 1965 .}

\item{[M.74]  Malgrange, B.  \textit{Int\'egrales Asymptotiques et Monodromie .} Ann. Scient. ENS 4-i\`eme s\'erie  7 (1974) p. 405-430 .}

\item{[Mi.68]  Milnor, J.  \textit{Singular Points of Complex Hypersurfaces .} Ann. of Math. Studies  61 (1968) Princeton  .}

\item{[P.86]  Palamodov, V.P.  \textit{Asymptotic expansions of integrals in complex and real regions .} Math. USSR 55 (1986) p.207-236  .}

\end{itemize}

\bigskip

\bigskip

\bigskip

\noindent Daniel Barlet , 

\noindent Universit\'e Henri Poincar\'e (Nancy I ) et Institut Universitaire de France,

\noindent Institut E.Cartan  UHP/CNRS/INRIA, UMR 7502 ,

\noindent Facult\'e des Sciences et Techniques, B.P. 239

\noindent 54506 Vandoeuvre-les-Nancy Cedex , France .

\end{document}